%
%

\magnification=1200

\font\titfont=cmr10 at 12 pt



\def\BBB{\bf}

\overfullrule=0in

\def\boxit#1{\hbox{\vrule
 \vtop{%
  \vbox{\hrule\kern 2pt %
     \hbox{\kern 2pt #1\kern 2pt}}%
   \kern 2pt \hrule }%
  \vrule}}
  \def\mathqed{  \vrule width5pt height5pt depth0pt}

  \def\harr#1#2{\ \smash{\mathop{\hbox to .3in{\rightarrowfill}}\limits^{\scriptstyle#1}_{\scriptstyle#2}}\ }

 \def\GG{{{\bf G} \!\!\!\! {\rm l}}\ }

\def\wee{\wedge\cdots\wedge}

\def\bra#1#2{\langle #1, #2\rangle}

\def\ss{\subset}

\def\half{\hbox{${1\over 2}$}}
\def\smfrac#1#2{\hbox{${#1\over #2}$}}

\def\dim{{\rm dim}}

\def\log{{\rm log}}
\def\Hess{{\rm Hess}}

\def\tr{{\rm tr}}
\def\max{{\rm max}}
\def\min{{\rm min}}
\def\span{{\rm span\,}}

\def\det{{\rm det}}
\def\End{{\rm End}}
\def\Sym{{\rm Sym}^2}

\def\Herm{{\rm Herm}}
\def\arr{\longrightarrow}

\def\Core{{\rm Core}}

\def\rn{\bbr^n}

\def\Int{{\rm Int}}

\def\Symn{{\Sym(\rn)}}

\def\Theorem#1{\medskip\noindent {\bf THEOREM \bf #1.}}
\def\Prop#1{\medskip\noindent {\bf Proposition #1.}}
\def\Cor#1{\medskip\noindent {\bf Corollary #1.}}
\def\Lemma#1{\medskip\noindent {\bf Lemma #1.}}
\def\Remark#1{\medskip\noindent {\bf Remark #1.}}
\def\Note#1{\medskip\noindent {\bf Note #1.}}
\def\Def#1{\medskip\noindent {\bf Definition #1.}}

\def\Ex#1{\medskip\noindent {\bf Example \bf    #1.}}

\def\pf{\medskip\noindent {\bf Proof.}\ }
\def\qed{\hfill  $\vrule width5pt height5pt depth0pt$}

\def\hk{\_{\rm l}\,}

   \def\cp{{\cal P}}

\def\ch{{\cal H}}   
   \def\cn{{\cal N}}

\def\cp{{\cal P}}
\def\cf{{\cal F}}

\def\vf{\varphi}

\def\wt{\widetilde}
\def\wh{\widehat}

\def\and{\qquad {\rm and} \qquad}
\def\arr{\longrightarrow}

\def\bbr{{\bf R}}\def\bbh{{\bf H}}
\def\bbc{{\bf C}}

\def\bbz{{\bf Z}}
\def\bbp{{\bf P}}

\def\a{\alpha}

\def\d{\delta}
\def\e{\epsilon}
\def\f{\phi}

\def\l{\lambda}
\def\o{\omega}

\def\s{\sigma}
\def\x{\xi}

\def\D{\Delta}
\def\L{\Lambda}

\def\O{\Omega}

\def\lag{Lagrangian}
\def\psh{plurisubharmonic }

\def\PH#1{\widehat {#1}}
\def\lloc{L^1_{\rm loc}}

\def\bo{\partial \Omega}

\def\PSH{{ \rm PSH}}

\def\lagpsh{Lagrangian plurisubharmonic}

\def\PSHl{\PSH_{\rm Lag}}

\def\Symn{\Sym(\rn)}
 
\def\USC{{\rm USC}}

\def\cpt{\wt{\cp}}
\def\ft{\wt F}
\def\ob{\overline{\O}}

 \def\LAG{{\rm LAG}}

 \def\ve{{\varepsilon}}
\def\cn{\bbc^n}
\def\sym{{\rm sym}}
\def\sk{{\rm skew}}

\def\cn{{\bbc^n}}
\def\Lag{{\rm Lag}}
\def\sk{{\rm skew}}
\def\sym{{\rm sym}}\def\cpL{{\cp(\LAG)}}
\def\id{{\rm I}}
\def\mlag{{\rm M}_{\Lag}}
\def\SS{{\bf S} \!\!\!\! /}
\def\bcl{{\rm C\!\!\! l \, l}_{2n}}

\def\PTCG{HL$_1$}
\def\DDD{HL$_2$}
\def\PinGGC{HL$_3$}
\def\DDR{HL$_4$}
\def\HLG{HL$_5$}
\def\GPandC{HL$_6$}
\def\HLGG{HL$_7$}
\def\SURVEY{HL$_8$}
\def\REST{HL$_9$}
\def\RSing{HL$_{10}$}
\def\ACM{HL$_{11}$}
\def\AE{HL$_{12}$}
\def\IDP{HL$_{13}$}
\def\EDGE{HL$_{14}$}
\def\CONT{HL$_{15}$}

\def\AA{1}
\def\BBB{2}
\def\BB{3}
\def\II{4}
\def\CC{5}
\def\DD{6}

\def\FF{7}
\def\HH{8}
\def\JJ{9}
\def\KK{10}
\def\LL{11}
\def\MM{12}
 
\ 
\vskip .3in

\centerline{{\titfont LAGRANGIAN POTENTIAL THEORY} AND A  }
\smallskip

\centerline{\titfont  LAGRANGIAN EQUATION OF MONGE-AMP\`ERE TYPE }
\bigskip

\centerline{\titfont F. Reese Harvey and H. Blaine Lawson, Jr.$^*$}
\vglue .9cm
\smallbreak\footnote{}{ $ {} \sp{ *}{\rm Partially}$  supported by
the N.S.F. }

\vskip .6in
\centerline{\bf ABSTRACT} \medskip
  \font\abstractfont=cmr10 at 10 pt

{{\parindent= .5 in\narrower\abstractfont 
The purpose of this paper is to establish a Lagrangian potential theory, analogous to 
the classical pluripotential theory, and to define and study a Lagrangian differential 
operator of Monge-Amp\`ere type.  This developement is new even in $\cn$.  However,  it applies
quite generally --  perhaps most importantly to symplectic manifolds equipped with 
a Gromov metric.

The Lagrange operator is an explicit polynomial on $\Sym(TX)$ whose principle branch
defines the space of Lag-harmonics.  Interestingly the operator depends only on the Laplacian and the
SKEW-Hermitian part of the Hessian.
The  Dirichlet problem for this operator is solved
in both the homogeneous and inhomogeneous cases.  It is also solved for each of the
other branches.

This paper also introduces  and systematically studies the notions of Lagrangian plurisubharmonic and harmonic
functions, and Lagrangian convexity.  An analogue of the Levi Problem is proved. 
In  $\cn$ there is another concept,  Lag-plurihamonics, which relate in several ways to the harmonics
on any domain.
Parallels of this Lagrangian potential theory with standard (complex) pluripotential theory are constantly  emphasized.

}}

\vfill\eject

\vskip .5in

\centerline{\bf TABLE OF CONTENTS} \bigskip

{{\parindent= .1in\narrower\abstractfont \noindent

\qquad \AA. Introduction.\smallskip

\qquad\BBB.  Four Fundamental Concepts.  \smallskip

\qquad \BB.  The Lagrangian Subequation and Lagrangian Hessian \smallskip

\qquad \II.   Lagrangian Pluriharmonic Functions  \smallskip

\qquad \CC.  A  Lagrangian Operator of Monge-Amp\`ere Type  \smallskip

\qquad \DD.  Two  Intrinsic Definitions of the Lagrangian  (Monge-Amp\`ere)  Operator \smallskip

\qquad \FF. Lagrangian Potential Theory on Gromov Manifolds   \smallskip

\qquad \HH.  The  Lagrangian Monge-Amp\`ere  Operator on Gromov Manifolds   \smallskip

\qquad \JJ.   Lagrangian Pseudoconvexity   \smallskip

\qquad \KK. Lagrangian Boundary Convexity   \smallskip

\qquad \LL. The Dirichlet Problem for the Lagrangian Monge-Amp\`ere Equation.  \smallskip

\qquad  \MM. Ellipticity of the Linearization.   

\vskip.3in
 \qquad Appendix A.  A More Detailed Presentation of the Lagrangian Subequation.
}}

\vfill\eject

\centerline{\bf \AA. Introduction.}
\medskip

The aim of this paper is to establish a Lagrangian potential theory, analogous to 
the classical pluripotential theory, and to define and study a Lagrangian differential 
operator of Monge-Amp\`ere type.  This developement is new even in $\cn$.  However,  it applies
quite generally --  perhaps most importantly to symplectic manifolds equipped with 
a Gromov compatible metric, which we shall call {\sl Gromov manifolds}.

On $\cn = (\bbr^{2n}, J)$  the Lagrangian operator  $\mlag(D^2 u)$ is a homogeneous polynomial in 
$A  = D^2 u  \in \Sym_\bbr(\bbr^{2n})$ of degree $N = 2^n$.  Its principal branch -- the closure of the 
component of $\{\mlag(A)\neq 0\}$ containing
the identity matrix -- has a fundamental geometric interpretation as the set of $A$'s satisfying
$$
\tr\left\{  A\bigr|_W   \right\} \ \geq\ 0 \qquad {\rm for\ all\  Lagrangian \ } \ n\ {\rm planes}\ \ W.
\eqno{(\AA.1)}
$$
Interestingly, this operator only depends on $\tr(A)$ and the skew-Hermitian part $A^\sk \equiv \half(A+JAJ)$ of $A$.
The value is unchanged if one adds to $A$ any Hermitian symmetric matrix $B=\half(B-JBJ)$ of trace 0.
There are several constructions of this operator, two of which have  intrinsic interpretations.
Either of these allows one to define the operator on any Gromov manifold $X$.  
One construction  uses the basic Spin$^c$-bundle
$\L^{0, *}(X)$.  The other comes from a certain derivation on the bundle $\L^n_\bbr \bbr^{2n}$.

One of the basic results here establishes  the existence and uniqueness of the solution
 to the Dirichlet problem for both the homogeneous and
the inhomogeneous Lagrangian Monge-Amp\`ere operator.  In $\cn$, the existence holds for  any smoothly
bounded domain $\O\ss\ss \cn$ such that $\bo$ satisfies a  
strict {\sl Lagrangian convexity condition} at each point.  One very geometric
version of this condition is that the trace of the second fundamental form on any tangential Lagrangian
$n$-plane is strictly inward-pointing.

The operator $\mlag(A)$ is actually G\aa rding hyperbolic with respect to the identity $I \in  \Sym_\bbr(\bbr^{2n})$,
which implies  that it has $N\equiv 2^n$ nested subequations corresponding to the increasing G\aa rding eigenvalues
$\L_1(A) \leq \L_2(A) \leq \cdots \leq \L_N(A)$.  The condition $\L_k \geq 0$ is the {\sl $k^{\rm th}$-branch} of the 
Lagrangian Monge-Amp\`ere operator. (The homogeneous equation for the operator above is the first branch.)
For each $k=1, ... , N$,   it is shown that  solutions to the homogenous Dirichlet problem for the equation 
$\L_k(D^2 u) =0$ are unique on all bounded domains $\O\ss\ss\cn$ as above.  
The appropriate boundary convexity becomes less stringent as $k$ decreases from $N$ to $N/2$.

Furthermore, all of these results concerning the Dirichlet problem carry over to any Gromov manifold
$X$. The statements are essentially the same.  For $\mlag$ one considers domains $\O\ss\ss X$ with smooth, strictly
Lagrangian-convex boundary.  The only additional hypothesis, which is global and always  true for $X = \cn$,
 is that there exist some smooth strictly Lagrangian plurisubharmonic function 
 defined on a neighborhood of $\ob$. (See the following paragraph 
for the terminology here.)

The original motivation for this study was to develop Lagrangian potential theory as another
tool in complex and symplectic analysis.  For this one can start with the notion of a {\sl smooth Lagrangian
plurisubharmonic} function defined on $\O\ss \cn$.  This is one whose second derivative (or  riemannian Hessian in the
general case) satisfies (\AA.1), i.e., that its trace on every Lagrangian $n$-plane at every point
is $\geq 0$. This notion of Lagrangian plurisubharmonic (or Lag-psh) can be carried over to the class of upper semi-continuous functions by using viscosity test functions.  The result  enjoys all the basic properties of classical plurisubharmonic functions.

One can show that Lag-psh functions are subharmonic, and so, in particular they are in $\lloc$.
There is also a Restriction Theorem which says that the restriction of a (viscosity) Lag-psh function
to any minimal Lagrangian submanifold $M$ is subharmonic on $M$ (or $\equiv -\infty$).
We note, by the way, the converse in $\cn$, that  if its restriction to
every affine Lagrangian plane is subharmonic, then it is (viscosity)  Lag-psh.  (This result, as in complex analysis, 
justifies the word ``plurisubharmonic''.)

We now return  to smooth functions and assume that $X$ is non-compact and connected.
Then given a compact subset $K\ss X$ we define the {\sl Lagrangian hull $\wh K$
of $K$} to be  those points $x\in X$ such that $u(x) \leq \sup_K u$ 
for all  Lag-psh  $u\in C^\infty(X)$.  Then $X$ is defined to be {\sl Lagrangian convex} if 
$K\ss\ss X \Rightarrow \wh K\ss\ss X$. This is meaningful on any non-compact Gromov manifold.           
 We prove that $X$ is Lagrangian convex iff $X$ admits a 
smooth proper exhaustion function which is Lag-psh.  
Moreover, if $X$ carries some smooth strictly Lag-psh function (for example, when $X=\cn$), then this exhaustion can be made  strict.
In this case $X$ has the homotopy-type of a 
complex of dimension $\leq 2n-2$. 
 
A  submanifold $M\ss X$ is defined to be {\sl Lag-free} (or free of lagrangian tangents)
 if it has no  tangential Lagrangian planes.
 (Thus, any manifold of dimension $< n$ is automatically free.)  If $M\ss X$ is  closed and free, then
 the domains $M_\e \equiv \{x : {\rm dist}(x, M)\leq \e\}$ admit strictly Lag-psh exhaustions
  for all $\e>0$ sufficiently small.

In Section 10 we discuss the notion of  strict  Lagrangian  boundary  convexity (referred to above).  The notion has several  equivalent definitions.  A local-to-global theorem is proved -- namely, if a domain $\O\ss X$ has 
a strictly Lagrangian convex boundary, then there exists a smooth exhaustion function which is strictly
Lag-psh at infinity. 
 (If, in addition, $X$ carries some strictly Lag-psh function, the exhaustion can be 
taken to be strictly  Lag-psh everywhere.)  This all can be viewed
as a (weaken) form of the Levi-problem in complex analysis.

Let us return now to the Lagrangian Monge-Amp\`ere equation.  Let 
$$
\cp(\LAG) \ \equiv\ \left \{A\in \Sym_\bbr(\cn) : \tr\left(A\bigr|_W    \right) \geq 0 \ \ \forall\, {\rm Lagrangian}\ \ W \right \},
$$
and note that a $C^2$-function $u$ is Lag-psh if $D^2u(x)\in \cp(\LAG)$ for all $x$  (and it is strictly Lag-psh if
$D^2u(x)\in \Int \cp(\LAG)$  for all $x$).   Now for a $C^2$-function $u$ it is natural to define $u$
   to be {\sl Lagrangian harmonic} if equality holds in the sense that
 $D^2u(x) \in \partial \cp(\LAG)$ for all $x$. Using the dual subequation 
 $\{\L_N\geq 0\}$ to   $\cp(\LAG) = \{\L_1\geq 0\}$, this notion can be carried over to 
continuous functions, and all of this makes sense on a Gromov manifold by replacing $D^2u$ with  $\Hess\, u$,
the riemannian Hessian of $u$.  These Lagrangian harmonic functions are exactly the  solutions to the principle branch
of the Lagrangian Monge-Amp\`ere equation.  We also see that these solutions, or Lag-harmonic functions,  are the {\sl maximal functions} in our Lagrangian pluripotential theory.

By the way, there is another notion, distinct from those of Lag-psh, dually Lag-psh and Lag-harmonic, 
which is the analogue of pluriharmonic in complex analysis.
 Namely, a function $u$ is {\sl Lagrangian pluriharmonic} if its restriction to each affine Lagrangian plane
 is classically harmonic.  
 (There is a notion of the Lagrangian Hessian, denoted $D_\LAG^2$, and $u$ is Lagrangian pluriharmonic if and only if
 $D_\LAG^2 u = 0$.)
 In Section \II \  the Lagrangian pluriharmonics  are characterized and then put to use by obtaining
a  classical (non-viscosity) interpretation of the notion of dually Lag-psh (see Theorem \II.6).
There is also an interesting result about their contact sets with Lag-harmonics (see Theorem \II.8).

Section 11 is concerning with the Dirichlet Problem for the Lagrangian  Monge-Amp\`ere operator
for domains in a Gromov manifold.  The results are valid for the full inhomogenous  case $\mlag(u) = \psi$
where  $\psi\in C(\ob)$ satisfies $\psi\geq 0$.  This result also carries over to certain operators associated
to the branches  (see Note \LL.5)
${\rm M}_{\Lag,k}(u) = \psi$,  again with $\psi\geq 0$.

In Section 12 we give a simple proof that the linearization of the operator at a smooth strictly $\LAG$-psh function
on a compact set, is uniformly elliptic.


\vskip .4in

\centerline {\bf \BBB. Four Fundamental Concepts.}   
\medskip

We first consider complex euclidean space $V\equiv \bbc^n$ with all of its standard structures: the complex structure
$J$, the euclidean inner product $\bra \cdot\cdot$, and the hermitian inner product 
$(\cdot, \cdot) \equiv \bra \cdot\cdot -i \o (\cdot, \cdot)$ where $\o(v,w) \equiv \bra{Jv}{w}$ is the standard
symplectic form.
Given a smooth function $u$, its hessian, or second derivative $A\equiv D^2_z u \in \Sym_\bbr(V)$
is a real symmetric bilinear form on $V$, or equivalently (by polarization)  a real quadratic form on $V$.
It can also  be considered as a linear map $A\equiv D_z^2 u\in\End_\bbr(V)$ which is symmetric.
All three of these natural isomorphisms  are used without mention throughout the paper,
and when an orthonormal basis is present, $A\equiv D_z^2u$ can also be considered a symmetric
$2n\times 2n$-matrix.

We emphasize the parallels with complex potential theory where the focus is on $\bbc$-plurisubharmonic
functions, that is, real-valued functions $u$ with
the property that the restriction of $u$ to complex lines is $\D$-subharmonic. Equivalently,
the restriction of the quadratic form $A\equiv D^2_zu$ to each complex line $W$ has non-negative 
trace, i.e., $\tr\left( A\bigr|_W\right)\geq0$ for all $W\in G_\bbc(1,\bbc^n) = \bbp_\bbc^{n-1}$, 
the Grassmannian of complex lines in $\bbc^n$.

Now our Lagrangian case can be presented in a manner analogous to the complex case. 
The subset $\LAG\ss G_\bbr(n,\cn)$
of Lagrangian $n$-planes in $\cn$ is defined by 
$$
W\in \LAG \quad\iff \quad JW= W^\perp  \quad\iff \quad \o\bigr|_W\ =\ 0.
\eqno{(\BBB.1)}
$$

\bigskip
\centerline{\bf Smooth Functions}
\medskip

Throughout this subsection $u$ is a smooth function defined on an open subset of $\cn$.
Note that for any real affine subspace $W$ of $\cn$ we have that
$
\tr\left(   D^2 u\bigr|_W  \right) \ =\ \D\left( u\bigr|_W\right)
$.

\Def{\BBB.1}  
We say that $u$ is  {\bf Lagrangian plurisubharmonic}
if at each point $x$ in its domain and for each affine Lagrangian $n$-plane  $W$ containing $x$, one has
$$
\tr\left(   D^2 u\bigr|_W  \right)   \ =\   \D\left( u\bigr|_W\right)\ \geq\ 0.
\eqno{(\BBB.2)}
$$
The constraint set on the second derivative $A=D^2_z u$ will be denoted by
$$
\cp(\LAG) \ \equiv \ \left \{A\in \Sym_\bbr(\cn) : \tr\left( A\bigr|_W\right)\ \geq 0 \ \ \forall\, W\in\LAG\right\}.
\eqno{(\BBB.3)}
$$

Note that this defines a convex cone in $\Sym_\bbr(\cn)$,
and that the definition can be reformulated by requiring 
$$
D^2_z u\in \cp(\LAG) \quad{\rm for\ all}\ \ z\ \ {\rm in\ the\ domain\ of}\ \ u.
\eqno{(\BBB.2)'}
$$
\medskip

The next concept  is central to the Lagrangian Dirichlet Problem
(and some possible notions of Lagrangian capacity, which we do not discuss). 
It is analogous to the notion of a {\sl maximal function}
or a solution of the complex Monge-Amp\`ere equation in complex potential theory (cf.\ [BT$_*$]).

\Def{\BBB.2}  We say that $u$  is {\bf Lagrangian  harmonic} if it is Lagrangian plurisubharmonic and 
 if at each point $z$ in its domain there exists an affine Lagrangian $n$-plane through $z$ such that
$$
\tr\left( D^2 u\bigr|_W\right) \ =\ \D \left(  u\bigr|_W\right) \ =\ 0.
\eqno{(\BBB.4)}
$$
In terms of the constraint set $\cpL$ this definition is equivalent to the requirement that
$$
D^2_z u \in \partial \cpL \ \ 
{\rm for\ all}\ \ z\ \ {\rm in\ the\ domain\ of}\ \ u.
\eqno{(\BBB.4)'}
$$
\medskip

In complex potential theory a function $u$ is said to be pluriharmonic if its restriction to each
complex line is harmonic.  Analogously we have the following concept.

\Def{\BBB.3}  We say that $u$ is {\bf Lagrangian pluriharmonic} if 
$$
\tr\left( D^2 u\bigr|_W\right) \ =\ \D \left(  u\bigr|_W\right) \ =\ 0  \ \ \ {\rm for \ all\ affine \ Lagrangian\ planes } \ W.
\eqno{(\BBB.5)}
$$
In terms of the constraint set $\cpL$ this definition is equivalent to the requirement that
$$
D^2_z u \in \cpL  \cap (-\cpL) \ \equiv \ E  
\qquad 
{\rm for\ all}\ \ z\ \ {\rm in\ the\ domain\ of}\ \ u.
$$
We note that $E$ is a vector subspace of the convex cone $\cp(\LAG )$, and it contains all
other vector subspaces  of $\cp(\LAG )$.
This set $E$ is called the {\bf edge of the convex cone} $\cp(\LAG )$.
\medskip

\bigskip
\centerline{\bf Upper Semi-Continuous Functions}
\medskip

Each of these three concepts can be extended from smooth function to 
the general class of upper semi-continuous functions by using viscosity test functions.
We begin with some generalities.
Given $X\ss\bbr^N$, let $\USC(X)$ denote the set of  upper semi-continuous
functions $u:X \to \bbr\cup\{-\infty\}$.

\Def{\BBB.4}
(a)  For  $u\in \USC(X^{\rm open})$ and  $x\in X$, a smooth function $\vf$, defined
on a  neighborhood of $x$, is called a {\bf (viscosity) test function} for $u$ at $x$ if
$$
u\ \leq\ \vf \quad{\rm near} \ \ x\ \ {\rm with\ equality\  at\ \ } x.
\eqno{(\BBB.6)}
$$
\medskip

\noindent
(b)  A closed subset $F\ss\Sym(\bbr^N)$ is called a {\bf subequation} if $F+\cp\ss F$.

\medskip
\noindent
(c)  Given  $u$ and $F$ as above, we say that  $u$  is {\bf $F$-subharmonic} if 
$$
D^2_x \vf\in F,
\eqno{(\BBB.7)}
$$
for each $x\in X$ and each test function $\vf$ for $u$ at $x$.

\medskip
\noindent
(d)  Given two subequations $F$ and $G$, the set 
$$
E \ \equiv\ F\cap (-G)\quad{\rm  is\ called\ a\ }
{\bf generalized\ equation}.
\eqno{(\BBB.8)}
$$

\medskip
\noindent
(e)  A continuous function $u$ is a  {\bf solution to $E$} if
$$
U \ \ {\rm is\ } F\,{\rm subharmonic\ and\ \ } -u \ \ {\rm is\ } G\,{\rm subharmonic}
\eqno{(\BBB.9)}
$$

\medskip
\noindent
(f)  Given a subequation $F$, the {\bf dual subequation} is defined to be
$$
\ft\ \equiv\ -(\sim \Int F) \ =  \   \sim(-\Int F)
\eqno{(\BBB.10)}
$$

\medskip
\noindent
(g)  A generalized equation of the form $\partial F = F\cap (-\ft)$ is called an {\bf equation}, and its
solutions are called {\bf $F$-harmonic}. Note that $u$ is a solution if and only if $u$ is $F$-subharmonic and
$-u$ is $\ft$-subharmonic.
\medskip

Note that, in particular,  we have  the dual subequation $\wt \cp(\LAG)$ of $\cp(\LAG)$.
It is given by
$$
\wt \cp(\LAG) \ =\ \left\{A\in \Sym_\bbr(\cn) : \exists \, W\in\LAG \ {\rm such\ that\ } \tr \left(A\bigr|_W\right)\geq0\right\}
\eqno{(\BBB.11)}
$$
Now we can define our fourth fundamental concept. 
First, a smooth function $u$ is {\bf dually Lagrangian plurisubharmonic} if $D^2_z u\in \cpt(\LAG)$
 at each point $z$ in its domain, that is, at each $z$ there exists $W\in \LAG$ such that
 $\tr \{D^2_z u\bigr|_W\geq0\}$

It is now easy to extend all four concepts for smooth functions  to upper semi-continuous functions on 
an open set $X\ss \cn$.

\Def{\BBB.5} Given $u\in \USC(X)$, we say that:
\medskip
\noindent
(1) $u$ is {\bf Lagrangian plurisubharmonic} if $u$ satisfies the subequation $\cpL$.

\medskip
\noindent
(2) $u$ is {\bf dually Lagrangian plurisubharmonic} if $u$  satisfies the subequation $\cpt(\LAG)$.
\medskip
\noindent
(3) $u$ is {\bf Lagrangian harmonic} if $u$ is a solution to the equation $\partial \cpL$.
Note that this holds if and only if $u$ is  Lag-psh and $-u$ is dually  Lag-psh.

\medskip
\noindent
(4) $u$ is {\bf Lagrangian pluriharmonic} if $u$ is a solution to the generalized equation 
$E\equiv  \cpL \cap (-\cpL)$, i.e., $u$ is Lag-psh and $-u$ is also Lag-psh.

\smallskip\noindent
The more succinct expressions are  respectively that $u$ is \medskip
\centerline{
(1) {\bf  Lag-psh} \quad
  (2) {\bf dually Lag-psh}\quad 
 (3) {\bf Lag-harmonic}\quad
  (4) {\bf Lag-pluriharmonic.}
}
\medskip

The terminology {\sl Lagrangian plurisubharmonic} is justified by the 
Restriction Theorem (see [\REST]).

\Theorem{\BBB.6}
{\sl
An upper semi-continuous function $u$ is Lag-psh if and only if its restriction
to each affine Lagrangian $n$-plane is $\D$-subharmonic (or $\equiv -\infty$).}
\medskip

This is directly analogous to the subequation $\cp$ for convex functions (restrict to affine lines)
and to the subequation $\cp(G_\bbc(1,\cn))$ for plurisubharmonic functions (restrict to affine complex lines).

  In Section \II \ we give a classically formulated (i.e., non-viscosity) characterization of dually Lag-psh functions, which complements  Theorem \BBB.6
  
  The following is  a useful fact.
  \Prop{\BBB.7}
$$
\cp(\LAG)\ \ss\ \D\ \equiv\ \{\tr A\ \geq\ 0\},
\eqno{(\BBB.12a)}
$$
or equivalently,
$$
{\rm Each \ Lag\!\!-\!\! psh \ function\ is\ classically\ subharmonic.}
\eqno{(\BBB.12b)}
$$
\pf
The $x$-axis $\rn$ and the $y$-axis $i\rn$ in $\cn$ are both Lagrangian $n$-planes.
They define the ``partial'' Laplacians 
$$
\D^x \ \equiv \  \left \{A : \tr \left( A \bigr|_{\rn} \right)\geq0  \right  \}
\and 
\D^y \ \equiv \  \left \{A : \tr \left( A \bigr|_{i\rn} \right)\geq0  \right  \}.
$$
Since $\cpL \ss \D^x\cap \D^y\ss \D$, (\BBB.12a) follows easily.\qed

\medskip

We conclude this section by pointing out that for any domain $\O\ss \cn$ the set, denoted
$\cp_{\Lag}(\O)$,  of all (u.s.c.) Lag-psh functions on $\O$,  has all the classical properties, namely:

\medskip
(1) \ \  $u\in \cp_{\Lag}(\O) \ \Rightarrow u\in\lloc$ (by (\BBB.12b)),

\medskip
(2) \ \  $u, v \in \cp_{\Lag}(\O) \ \Rightarrow \max\{u,v\}\in \cp_{\Lag}(\O)$,

\medskip
(3)\ \ $\cp_{\Lag}(\O)$ is closed under uniform limits and decreasing limits, 

\medskip
(4) \ \ if $\cf \ss \cp_{\Lag}(\O)$ is any family locally bounded above, and if 

\qquad \  $U(x)\equiv \sup\{u(x) : u\in \cf\}$ for $x\in\O$, 

\qquad \   then the  u.s.c. regularization $U^*\in \cp_{\Lag}(\O)$, and

\medskip
(5) \ \ a function $u\in C^2(\O) \cap \cp_{\Lag}(\O))$ is Lag-psh in the sense of Definition \BBB.1.

\medskip
The properties  (2) -- (5) have been well known for a long time for any subequation. Proofs can be found, for example, in Appendix B in [\DDR]. Finally note the {\sl convex composition property}:

\medskip
(6) \ \  If $u\in \cp_{\Lag}(\O)$ and if $\chi$ is a convex, increasing $\bbr$-valued function on  Image$(u)$,

\qquad \  
then $\chi\circ u\in \cp_{\Lag}(\O)$.

\medskip\noindent
This follows since $\cp(\LAG)$ is a convex cone (see Fact (2) in \S 6 of  [\PinGGC]).

\vfill\eject

\centerline {\bf \BB. The Lagrangian Subequation and Lagrangian Hessian.}   
\medskip

A more detailed algebraic (and self-contained) discussion of this subequation 
is presented in Appendix A. Here we 
again  emphasize the parallels with complex potential theory, where 
the pertinent  subequation is defined by the subset $\bbp_\bbc^{n-1} \ss G_\bbr(2, \cn)$
of complex lines in $\cn$ --
namely, the subequation:
$$
\cp(\bbp_\bbc^{n-1})\ =\ \left \{A\in\Sym_\bbr(\cn) : \tr\left( A\bigr|_W\right)\geq0\ \ 
\forall \ W\in \bbp_\bbc^{n-1}\right\}.
$$
  Note that for any unit vector $e\in W \in \bbp_\bbc^{n-1}$,
$$
\tr\left( A\bigr|_W\right) \ =\ \bra {Ae} e + \bra {AJe}{Je} \ =\ \bra{(A-JAJ) e}e.
\eqno{(\BB.1)}
$$
Consequently, the subspace ${\rm Herm}^{\rm skew}(\cn)\ss \Sym_\bbr(\cn)$ of real symmetric maps  $A$ 
which anti-commute with $J$ is unimportant here since, given $A\in \Sym_\bbr(\cn)$, 
$$
\tr\left (A\big|_W\right) \ =\ 0 \quad
\forall \, W\in \bbp_\bbc^{n-1}
\qquad\iff\qquad
A \in {\rm Herm}^{\rm skew}(\cn)
\eqno{(\BB.2)}
$$
So in  classical pluripotential theory it is best to focus on its orthogonal complement in 
$\Sym_\bbr(\cn)$, namely the space  ${\rm Herm}^{\rm sym}(\cn)$ of real symmetric maps
on $\cn$ which commute with $J$.
The opposite is true for $\cp(\LAG)$, defined by (\BBB.3), where $\bbp_\bbc^{n-1}$ is replaced by $\LAG$.

Since this algebra will be used here, we summarize as follows.
The space $\Sym_\bbr(\cn)$ decomposes as an orthogonal direct sum
(with respect to $\bra AB = \tr(AB)$) as
$$
\Sym_\bbr(\cn)\ =\ {\rm Herm}^{\rm sym}(\cn) \oplus {\rm Herm}^{\rm skew}(\cn).
\eqno{(\BB.3)}
$$
Note that the dimensions are $n(2n+1) = n^2 + (n^2+n)$.  The orthogonal projections are
defined by 
$$
A^{\rm sym} \ \equiv\ \half(A-JAJ)
\and
A^{\rm skew} \ \equiv\ \half(A+JAJ)
\eqno{(\BB.4)}
$$
for all $A\in \Sym_\bbr(\cn)$.
The projection of the hessian onto the hermitian symmetric part:
$$
\left( D^2 u\right)^{\rm sym} \in {\rm Herm}^{\rm sym}(\cn)
\eqno{(\BB.5)}
$$
is called the {\sl   complex part} of the hessian of $u$,
or just the {\sl complex hessian} in pluripotential theory.
This is the important part
of the second derivative for complex analysis, and the subequation $\cp(\bbp_\bbc^{n-1})$ 
can be defined by the constraint $(D^2 u)^\sym\geq0$.  Also note that $\Herm^\sk(\cn)
= \cp(\bbp_\bbc^{n-1})  \cap (-\cp(\bbp_\bbc^{n-1}))$ is the edge of the convex cone $\cp(\bbp_\bbc^{n-1})$.

\medskip
\noindent
{\bf Note (Complex Notation).}
Besides the four ways of looking at an element $A\in\Sym_\bbr(\cn)$
(discussed at the the beginning of Section 2), there are additional identifications 
that come into play because of the complex structure. First, ${\rm Herm}^{\rm sym}(\cn)$ can be 
identified with the space  of complex linear maps $A\in \End_\bbc(\cn)$ which are self-adjoint,
as well as the space of complex matrices $A\i{\rm M}_n(\bbc)$ with $A= {\overline A}^t$.
In coordinates 
$$
(D^2 u)^\sym \ \cong\  \left( \partial^2 u \over \partial z_i \partial \bar z_j  \right)
\qquad {\rm is\ the\ complex\ hessian}.
\eqno{(\BB.6)}
$$
 The second space ${\rm Herm}^{\rm skew}(\cn)$ is isomorphic 
as a real vector space to $\Sym_\bbc(\cn)$, the space of symmetric $n\times n$ complex matrices,
or equivalently, the space of pure complex quadratic forms on $\cn$.  
The inverse of this isomorphism is given by sending a complex symmetric form $B(w) \equiv \sum_{ij} b_{ij} w_i w_j$
to $A\in \Sym_\bbr(\cn)$ defined by $A(w) = {\rm Re} \, B(w)$.  This gives the coordinate expression
$$
\left( D_z^2 u  \right)^\sk(w) \ 
=\  {\rm Re}\left( \sum_{i,j=1}^n {\partial^2 u\over \partial z_i \partial z_j} (z)w_i w_j    \right) \ =\  {\rm Re} \, B(w)
\eqno{(\BB.7a)}
$$
which shall be abbreviated to:
$$
\left( D_z^2 u  \right)^\sk  \ =\ {\bf Re}\left( {\partial^2 u\over \partial z_i \partial z_j}  \right).
\eqno{(\BB.7b)}
$$

\noindent
{\bf  Cautionary Note.}
Keep in mind that (\BB.7b)  denotes the real part of the quadratic form $ {\rm Re} \, B(w)$ in (\BB.7a), 
and {\bf not} the real part of the matrix $\left( {\partial^2 u\over \partial z_i \partial z_j}  \right)$.
Note also that 
$$
 {\rm Re} \, B(w) \ =\ 0 \qquad\iff\qquad
 B(w) \ =\ 0.
$$
The implication $\Rightarrow$ follows since ${\rm Re} \, B(\sqrt{i}w) = {\rm Re} \, iB(w) = -{\rm Im}\, B(w)$.

\medskip

Now the Lagrangian case can be presented in a manner analogous to the complex case. 
We first describe the Lagrangian analogue of (\BB.2). For this a standard canonical form is
needed, which for future reference we state explicitly.

\Lemma {\BB.1}
{\sl
For each $B\in \Sym_\bbr(\cn)$ which is skew hermitian (i.e., $BJ=-JB$)
there exists a unitary basis $e_1,...,e_n$ of $\cn$
and numbers $\l_j\geq0, \ j=1,...,n$ so that $B$ takes the canonical form
$$
B \ \equiv\ \l_1\left( P_{e_1} - P_{Je_1}  \right) + \cdots +  \l_n\left( P_{e_n} - P_{Je_n}  \right),
\eqno{(\BB.8)}
$$
where $P_v$ denotes orthogonal projection onto the line generated by $v$.
}

\pf
Since $B$ is skew hermitian, if $e$
is an eigenvector with eigenvalue $\l$, then $Je$ is an eigenvector with eigenvalue $-\l$.
\qed
\medskip


As with any unitary basis,  note that 
$$
W\ \equiv\ \span\{e_1, ... , e_n\} 
\and
JW\ \equiv\ \span\{Je_1, ... , Je_n\} 
\qquad {\rm are\ Lagrangian}
\eqno{(\BB.9)}
$$
Using this fact along with (\BB.8), we will  answer the natural question: which real symmetric forms
$A$ have the property that their restriction to every Lagrangian $n$-plane has trace zero?

Let ${\rm Herm}^{\rm sym}_0(\cn)$ denote the space of {\rm traceless} hermitian symmetric forms.
Note that given  $A\in \Sym_\bbr(\bbc^n)$, 
$$
A\in {\rm Herm}^{\rm sym}_0(\cn)
\qquad\iff\qquad
 A^\sk =0 \ \ \  {\rm and}\ \ \ \ \tr A=0.
$$
In particular, the decomposition of $\Sym_\bbr(\bbc^n)$ in (\BB.3) further decomposes as 
$$
\Herm^\sym(\cn) = [I]\oplus  {\rm Herm}^{\rm sym}_0(\cn).
\eqno{(\BB.3b)}
$$

\Lemma{\BB.2}
{\sl  
Suppose $A\in \Sym_\bbr(\cn)$. Then 
$$\tr\left(   A \bigr|_W  \right) = 0 \ \ \ \forall\, W\in\LAG 
\qquad\iff \qquad A\in {\rm Herm}^{\rm sym}_0(\cn).
$$
Said differently, $E\equiv \cpL \cap (-\cpL) = \Herm^\sym_0(\cn)$
is the generalized equation defining Lagrangian pluriharmonic functions,
cf. (\BB.5).}

\pf $(\Rightarrow)$ If $W = \span\{e_1, ... , e_n\}$ is the Lagrangian plane in Lemma \BB.1 for 
$B=A^\sk$,   then since each $\l_j\geq0, j=1,..., n$,
$$
\l_1 +\cdots + \l_n \ =\ \tr\left(A\bigr|_W\right)\ =\ 0
\qquad\Rightarrow\qquad A^\sk\ =\ 0.
$$
Moreover, for any $W\in\LAG$
$$
\tr A\ =\  \tr\left(A\bigr|_W\right) + \tr\left(A\bigr|_{JW}\right),
\eqno{(\BB.10)}
$$
proving that $\tr A=0$.

 $(\Leftarrow)$ Conversely, if $A\in {\rm Herm}^{\rm sym}_0(\cn)$, then
 $\bra {Ae} e = \bra {AJe}{Je}$, so that  $\tr\left(A\bigr|_W\right) =  \tr\left(A\bigr|_{JW}\right)$
 for all $W\in\LAG$. This proves, by (\BB.10) that $\tr\left(A\bigr|_W\right) =  0 \ \forall W\in\LAG$.\qed

\vskip .3in

\centerline{\bf
The Lagrangian Hessian
}

\medskip

Because of Lemma \BB.2, the orthogonal complement of 
${\rm Herm}^{\rm sym}_0(\cn)$ is now our focus. This space is exactly
$$
{\rm Herm}^{\rm sym}_0(\cn)^\perp \ =\  [{\rm I}] \oplus {\rm Herm}^{\rm skew}(\cn)
\eqno{(\BB.11)}
$$ 
where $ [{\rm I}]  \equiv \bbr \cdot {\rm I}$ denotes the line through the identity.
The projection of $A\in \Sym_\bbr(\cn)$ onto this space, which will be denoted $A^{\rm Lag} \equiv \pi(A)$
is given by
$$
A^{\rm Lag} \ \equiv\ {(\tr A)\over 2n} \cdot {\rm I}+ A^{\rm skew}
\eqno{(\BB.12)}
$$  
and is called the {\bf Lagrangian part of $A$}.

\Cor{\BB.3} {\sl
For any $A\in \Sym_\bbr(\cn)$ and $W\in\LAG$ one has that}
$$
\tr\left(A\bigr|_W\right) \ =\ \tr\left(A^\Lag\bigr|_W\right) \ =\ \half  \tr(A) +  \tr\left(A^{\rm skew}\bigr|_W\right)
\eqno{(\BB.13)}
$$ 
\Def{\BB.4. (The Lagrangian Hessian)} Suppose $u$ is a smooth function defined on an open
subset of $\rn$.  The  {\bf  Lagrangian hessian} of $u$ is the Lagrangian part of the second derivative
of $u$, that is:
$$
\Hess^\Lag (u) \ =\ (D^2 u)^\Lag \ =\ {\D u\over 2n}\cdot {\rm I} + \left(D^2 u\right)^{\rm skew}
\eqno{(\BB.14)}
$$ 
where $(D^2 u)^{\rm skew} \equiv \half (D^2 u +JD^2 u J)$.

\medskip
Note that by Lemma \BB.2, $u$ is lagrangian pluriharmonic  if and only if $\Hess^\Lag (u) \equiv 0$.
In particular,  the notion of Lagrangian plurisubharmonic only depends on the 
Lagrangian hessian.  Also note that just as the complex hessian can be expressed
in $z, \bar z$ coordinates as $(D^2u)^{\rm sym} \cong \left( {\partial^2 u\over \partial z_i \partial \bar z_j}   \right)$,
the Lagrangian hessian at a point $z$ can be written as the quadratic form
$$
\left(D^2_z\right)^\Lag \
 \cong \ \left({\D u\over 2n} \right) {\rm I} + {\bf Re}\left( {\partial^2 u\over \partial z_i \partial z_j}  \right)
\eqno{(\BB.15)}
$$ 
From (\BB.13) we see that
$$
\inf_{W\in\LAG}  \tr\left(A\bigr|_W\right)\ =\ \inf_{W\in\LAG}  \tr\left(A^\Lag\bigr|_W\right)
\eqno{(\BB.16)}
$$ 
which of course applies to $A=D^2 u$.

\vskip .3in

\centerline{\bf The Canonical Form for the Lagrangian Hessian}

\medskip

As in the complex case, once the key part of the second derivative has been identified,
it can be put in canonical form.

\Cor{\BB.5}
{\sl
For each $H\in   [{\rm I}] \, \oplus \,{\rm Herm}^{\rm skew}(\cn)$, there exists a unitary basis
$e_1, Je_1, ... e_n, Je_n$ for $\cn$ consisting of eigenvectors of $H$ with corresponding eigenvalues
$$
{\tr H\over 2n} +\l_1, \ {\tr H\over 2n} -\l_1, \ ... \ {\tr H\over 2n} +\l_n, \ {\tr H\over 2n} -\l_n.
\eqno{(\BB.17)}
$$ 
where  $\l_1\geq \l_2\geq\cdots\geq \l_n\geq0$ are the eigenvalues of $H^\sk$.}

\pf Note that $H=  {\tr H\over 2n} {\rm I} +H^\sk$ and $H^\sk$ has the canonical form (\BB.8).\qed

\medskip

Perhaps not surprisingly we have the following.

\Theorem{\BB.6}
{\sl
For each $A \in \Sym_\bbr(\cn)$
$$
\inf_{W\in \LAG} \tr\left( A\bigr|_W  \right) \ =\ {\tr A\over 2} - (\l_1 + \l_2+\cdots +\l_n)
\eqno{(\BB.18)}
$$ 
where $\l_1\geq\l_2\geq\cdots\geq\l_n\geq0$ are the non-negative eigenvalues of $A^\sk$.
}
\medskip

We shall give an indirect proof with a lemma which will be useful later.

First, for any $A\in \Sym_\bbr(\cn)$, we make the abbreviations: $\mu \equiv \half \tr A$,
$B \equiv A^\sk$, and $H\equiv A^\Lag \equiv {\mu\over n} I + B$ =  the Lagrangian part of $A$.
In the following choose the eigenstructure given by Corollary \BB.5
for $H\equiv A^\Lag$.

\Lemma{\BB.7}
{\sl
Given $H\equiv {\mu\over n}{\rm I} +B$ with $B\in {\rm Herm}^\sk(\cn)$,
extend $H$ as a derivation $D_H$ on $\L^n_\bbr\cn$.
Then the eigenvectors of $D_H$ are the ${2n}\choose n$ real axis $n$-planes in
$\bbr^{2n} =   \cn$ (in terms of the $B$-eigenstructure).  Each such axis $n$-plane $\x$ has a unique
decomposition $\x=\a\wedge \eta$ where $\a$ is an axis complex $p$-plane 
and $\eta$ is an isotropic axis $q$-plane (with $2p+q=n$).
The corresponding eigenvalue for $D_B$ is $\tr(  B\bigr|_{\span \eta} )$,
and for $D_H$ it is $ \mu+\tr(  B\bigr|_{\span \eta} )$.
}

\pf
First note that for each of the complex $p$-planes $\a$ we have $D_B\a=0$ since $D_B(e_i\wedge Je_i)=0$.
Any isotropic axis $q$-plane $\eta$ is of the form
$$
\eta\ =\ \{e_{i_1}\ {\rm or}\ J e_{i_1}\}\  \wedge \  \cdots \ \wedge\ \{e_{i_q}\ {\rm or}\ J e_{i_q}\}
\eqno{(\BB.19)}
$$ 
for some increasing multi-index $I=(i_1,...,i_q)$ of length $q$. With $+$ corresponding to the choice of 
$e_{i_j}$ and $-$ corresponding  to the choice of $Je_{i_j}$, we have
$$
D_B \eta \ =\ (\pm\l_{i_1} \pm \cdots\pm \l_{i_q})\eta
\eqno{(\BB.20)}
$$ 
Therefore the eigenvalues of $D_B$ are 
$$
\pm\l_{i_1} \pm \cdots\pm \l_{i_q} \qquad{\rm for \ all} \ \ I=(i_1,...,i_q),\ \ 0\leq q\leq n.
\eqno{(\BB.21)}
$$ 
This proves that the eigenvalues of $D_H$ are
$$
\mu \pm\l_{i_1} \pm \cdots\pm \l_{i_q} \qquad{\rm for} \ \ 0\leq |I|= q\leq n.
\eqno{(\BB.22)}
$$ 
since $D_{{\mu\over n}{\rm I}}\x =\mu\x$ for all axis $n$-planes $\x$.\qed\medskip

\noindent
{\bf Proof of  Theorem \BB.6.}   It is now obvious that the minimum eigenvalue
of $D_H$ equals 
$$
\l_{\rm min} \ =\ \mu-(\l_1+\cdots+\l_n)
\eqno{(\BB.23)}
$$ 
and the corresponding eigenvector is $\x= Je_1 \wee Je_n$.
Of course
$$
\l_{\rm min}(D_H) \ =\ \inf_L \tr\left( D_H\bigr|_L\right)
\eqno{(\BB.24)}
$$ 
where the inf is taken over all lines $L$ in $\L^n_\bbr \cn$. 
Finally note that each $W\in \LAG$ determines the line $L(W)$ in $\L^n_\bbr \cn$
through the simple $n$-form obtained by wedging a basis for $W$,  
and that
$$
\tr\left(D_H\bigr|_{L(W)}  \right)\ =\ \tr\left( H\bigr|_W  \right).
\eqno{(\BB.25)}
$$ 
(In fact (\BB.25) is true with $H$ replaced by any $A\in \Sym_\bbr(\cn)$ and $W$ any subspace.)
Therefore,
$$
\l_{\rm min}(D_H) \ =\  \inf_{W\in\LAG}\tr\left( H\bigr|_W  \right).  \qquad\mathqed
\eqno{(\BB.26)}
$$ 
\medskip

Combining (\BB.13), (\BB.26) and (\BBB.3) in  Definition \BBB.1  yields

\Theorem {\BB.8}
{\sl
With $H \equiv A^\Lag$, the Lagrangian part of $A\in\Sym_\bbr(\cn)$, one has
$$
A\in \cp(\LAG) \qquad\iff\qquad D_H\ \geq\ 0.
$$
In terms of smooth functions, $u$ is Lagrangian plurisubharmonic if and only if  at each
point $z$ in its domain, }
$$
D_H \ \geq\ 0\qquad{\rm where}\ \ H\ \equiv\ \left(D^2_z u \right)^\Lag
\quad{\rm is\ the\ Lagrangian\ Hessian \ of  \ } u.
$$

\Remark{\BB.9}
The operator  $f(D^2 u)$ defined by 
$$
f(A)\ =\ {\tr A\over 2} -\l_1^+- \cdots - \l_n^+,
\eqno{(\BB.27)}
$$ 
where  $\l_1^+, ..., \l_n^+ \geq0$ are the non-negative eigenvalues of $A^\sk$,
defines the subequation $\cpL$ by $f(A)\geq0$, and the equation $\partial \cp(\LAG)$ by 
also requiring that $f(A)=0$.
As discussed in [\HLG,  \S 3]  (see also [\DDR,  Rmk.\ 14.11]), each  subequation $G\ss\Symn$ is defined
as $G=\{A : g(A)\geq 0\}$
 for a unique operator $g:\Symn \to \bbr$ with the property that  $g(A+tI) = g(A) +t$, which we call the canonical operator for
the subquation $G$.  Since $g\equiv {1\over n} f$ with $f$ defined by (\BB.27) has this property,
the operator ${1\over n}f$ is the canonical operator for $G\equiv \cp(\LAG)$.

\vskip .3in

\centerline{ \bf  \II.  Lagrangian Pluriharmonic Functions.}   
\medskip

In this section we explore the analogues of the pluriharmonic functions in pluripotential theory
in our context of Lagrangian geometry.  
(In a more general context both are examples of  {\sl ``edge functions''} -- see Remark \II.4 below.)
Unlike their cousins,  Lagrangian pluriharmonics
do not exist on general manifolds.  Even in euclidean space our first result  (Theorem \II.1)
explains the limited nature of these functions.  However, in spite of this result, the 
Lagrangian pluriharmonics do play an important role by providing an equivalent way 
of defining the dually Lag-psh functions (Theorem \II.6),
which is formulated in classical language.

Recall that by Definition \BBB.6 a function $h\in C(\O)$, $\O^{\rm open}\ss\cn$,
is {\bf Lagrangian pluriharmonic} if both $h$ and $-h$ are Lag-psh.
By (\BBB.14), this implies that $\D h=0$ and hence $h$ is real analytic.
Much more is true.

\Theorem{\II.1}
{\sl
The space of Lag-pluriharmonic functions on a connected open subset $\O$
of $\cn$ consists of the traceless hermitian  degree 2 polynomials:
$$
h(z) \ =\ c+ \sum_{k=1}^n (b_k z_k +\overline{b}_k\overline{z}_k ) +{1\over 2}\sum_{i,j=1}^n a_{ij} z_i \overline{z}_j
\eqno{(\II.1)}
$$
where $c\in \bbc, b\in\cn$, and $A=(a_{ij})$ satisfies $A={\overline A}^t$ (see the first Note in \S3).
}
\pf 
The elementary proof is preceded by some elementary observations.  
One can also refer to a Lag-pluriharmonic as a solution to the edge equation
$$
D^2_z h\ \in\ E \ \equiv\ \cpL  \cap  (-\cpL)
\eqno{(\II.2)}
$$
Recall from Section \BB \ that $E = \Herm_0^\sym(\cn)$, and since $E^\perp = [\id]\oplus \Herm^\sk(\cn)$,
it is defined by the vanishing of the Lagrangian hessian $(D^2_z h)^\Lag$.
Since $\D h=0$, by (\BB.15) this is equivalent to 
$$
(D^2_z h)^\sk\ \cong \ {\bf Re}\left(  {\partial^2 h  \over \partial z_i \partial z_j}   \right)
\ =\ 0,
$$
which is equivalent to  
$
\left(  {\partial^2 h  \over \partial z_i \partial   z_j }  \right) =0
$
(see the Cautionary Note after (\BB.7)), 
This reduces the proof of Theorem \II.1 to the following lemma.

\Lemma{\II.2} 
{\sl
If $h$ is a smooth function satisfying
$$
(D^2_z h)^\sk\ \cong \   {\bf Re}\left(  {\partial^2 h  \over \partial z_i \partial  z_j}   \right)
\ =\ 0, \quad{\rm or\ equivalently\ \ }  {\partial^2 h  \over \partial z_i \partial  z_j}  \ =\ 0 \ \ \forall\, i,j,
$$
near a point $z$, then (\II.1) holds in a neighborhood of $z$.}


\pf 
As noted above 
$
  {\partial^2 h \over \partial z_i \partial  z_j } =0
$
for all $i,j$ is equivalent to $(D^2_z h)^\sk=0$.  Next note that 
$$
\eqalign
{
4   {\partial^2 h \over \partial z_i \partial z_j } 
&= \ \left( {\partial\over \partial x_i}  -i  {\partial\over \partial y_i}  \right) 
\left( {\partial h \over \partial x_j}  -i  {\partial h \over \partial y_j}   \right)     \cr
&=\  \left( {\partial^2 h\over \partial x_i  \partial x_j}  - {\partial^2 h\over \partial y_i  \partial y_j}  \right)
 -i  \left( {\partial^2 h\over \partial x_i  \partial y_j}  - {\partial^2 h\over \partial y_i  \partial x_j}  \right).
}
$$
Therefore,
$$
{\partial^3 h \over \   \partial x_i  \partial x_j  \partial x_k  }  \ =\ 
{\partial^3 h \over \   \partial y_i  \partial y_j  \partial x_k  }  \ =\ 
- {\partial^3 h \over \   \partial y_i  \partial x_j  \partial y_k  }  \ =\ 
- {\partial^3 h \over \   \partial x_i  \partial x_j  \partial x_k  }  
$$
all vanish.  Similarly,
$$
{\partial^3 h \over \   \partial y_i  \partial y_j  \partial y_k  }  \ =\ 
{\partial^3 h \over \   \partial x_i  \partial x_j  \partial y_k  }  \ =\ 
- {\partial^3 h \over \   \partial x_i  \partial y_j  \partial x_k  }  \ =\ 
- {\partial^3 h \over \   \partial y_i  \partial y_j  \partial y_k  }  
$$
all must vanish.  Thus all the third partial derivatives  vanish, proving
that $h$ is of degree 2.  Since the component $(D^2_z h)^\sk\ =0$,
we have $D^2 h\in \Herm^\sym(\cn)$, which proves (\II.1)
since $\D h=0$.\qed

\Remark{\II.3}
As noted above $E = \cpL\cap(-\cpL) \ =\ \Herm_0^\sym(\cn)$
is a generalized equation (see (\BBB.11)).  By contrast
$\Herm^\sym(\cn)$ is {\sl not} a generalized equation since the smallest
subequation $F$ containing $\Herm^\sym(\cn)$ is $F = \Herm^\sym(\cn)+\cp$,
which equals all of $\Sym_\bbr(\cn)$, because $\cp\cap  \Herm^\sk(\cn) = \{0\}$.

\Remark{\II.4.  (Edge Functions)}
For any convex cone subequation $\cp^+$, the {\sl edge}  is defined to be 
$E \equiv \cp^+ \cap (-\cp^+)$.  (See Appendix A for more details.)
In general, one can define, using  viscosity test functions, solutions of this
(generalized) equation.  It is natural to refer to such functions as {\sl edge functions},
or in the geometric cases (i.e., $\cp^+ \equiv \cp(\GG)$ with $\GG$ a closed subset 
of a Grassmannian) as {\sl $\GG$-pluriharmonic functions}.
In the four cases: $\cp^+ = \cp(G_\bbr(1, \rn)), \cp(G_\bbc(1,\cn)),  \cp(G_\bbh(1,\bbh^n))$ and 
$\cp(\LAG)$, the edges are: $E = \{0\}, \Herm_\bbc^\sk(\cn), \Herm_\bbh^\sym(\cn)^\perp$, and 
$\Herm_0^\sym$ respectively.   
The edge functions, or pluriharmonics, are:   affine functions  for $\cp(G_\bbr(1, \rn))$.
 Re$\{f(z)\}$ with $f$ holomorphic for $\cp(G_\bbc(1,\cn))$, and
in the case $\cp(\LAG)$, the traceless degree 2 polynomials (Theorem \II.1).
Thus the Lagrangian case is more like the convex case than the complex case in that the space of edge
functions is finite dimensional, whereas the space of real parts of holomorphic functions is infinite dimensional.

\medskip

The Lag-pluriharmonic functions are closely related to the dually Lag-psh functions
(i.e., the subharmonics for the dual subequation $\wt \cp (\LAG)$).

\Def{\II.5}  An upper semi-continuous function $u$ is ``sub'' the traceless hermitian 
degree-2 polynomials on an open subset $X\ss\cn$
if for all domains $\O\ss\ss X$ and all traceless hermitian 
degree-2 polynomials $h$,
$$
u\ \leq\ h\quad{\rm on}\ \ \bo 
\qquad\Rightarrow\qquad
u\ \leq\ h\quad{\rm on}\ \ \ob.
\eqno{(\II.4)}
$$

\Theorem{\II.6} {\sl
A function $u$ is dually \Lag-psh \ \ $\iff$\ \ $u$ is ``sub'' the  traceless hermitian 
degree-2 polynomials
 \ \ $\iff$\ \ $u$ is locally  ``sub'' the  traceless hermitian 
degree-2 polynomials.
}

\noindent
\pf Suppose $u$ is $\cpt(\LAG)$-subharmonic on $X$ and $h$ is a traceless hermitian
degree-2 polynomial.  Since $-h$ is $\cp(\LAG)$-subharmonic, (\II.4) follows from comparison
(see Thm. 6.2 in  [\AE]).

The proof of the converse illustrates the importance of Proposition A.2 in the Appendix.
Suppose now that $u$ is not $\cpt(\LAG)$-psh on $X$.
Then (see Lemma 2.4 in  [\DDR]) there exists $z_0\in X$, a quadratic polynomial $\vf$, and $\a>0$
such that
$$
u(z) \ \leq \ \vf(z) -\a|z-z_0|^2\quad {\rm near}\ \ z_0\ \ {\rm with\ equality\ at\ \ } z_0,
\eqno{(\II.5a)}
$$
but
$$
D^2_{z_0} \vf \notin \cpt(\LAG), \ \ \  {\rm i.e., }\ \ -D^2_{z_0}\vf \in   \Int\cpL.
\eqno{(\II.5b)}
$$
By Proposition A.2(4) we have $\Int \cpL = \Int \cp+\Herm_0^\sym(\cn)$. Thus
$$
 -D^2_{z_0}\vf \ =\ P+B  \ \ \ {\rm with}\ \ \ P>0\ \ \ {\rm and}\ \ \ B\in \Herm_0^\sym(\cn).
\eqno{(\II.5c)}
$$
The  traceless hermitian degree-2 polynomial $h$ satisfies
$$
\eqalign
{
h(z) \  &\equiv\  \vf(z_0) + \bra {D_{z_0} \vf}{z-z_0}  -  \half\bra {B(z-z_0)}{z-z_0}   \cr
&= \ \vf(z) -\half \bra  {D_{z_0}^2 \vf}{z-z_0}  -  \half\bra {B(z-z_0)}{z-z_0}   \cr
&=\  \vf(z) + \half \bra {P(z-z_0)}{z-z_0}.
}
$$
Therefore, (\II.5) implies 
$$
u(z) \ \leq\ h(z) - \a|z-z_0|^2
\eqno{(\II.6)}
$$
near $z_0$ with equality at $z_0$.  This implies that  $u$ is not sub the function $h$ 
on any small ball about $z_0$.  Hence, $u$ is not locally ``sub'' the traceless hermitian degree-2 polynomials.
\qed

\Remark{\II.7}  Theorem \II.6 has analogues for the subequations $\cp$ and its complex analogue.
For the first case it says that $u$ is $\cpt$-subharmonic if and only if it is ``sub'' the affine functions
(see [\DDD]). For the complex Monge-Amp\`ere subequation $\cp_\bbc \equiv \{A : A-JAJ\geq0\}$ on $\cn$, it says that 
$u$ is $\cpt_\bbc$-subharmonic if and only if $u$ is ``sub'' the pluriharmonic functions (see Prop.  5.14 in [\REST]).
\medskip

Fix an open set $X\ss \cn$ and a compact subset $K\ss X$.  Then the
Lagrangian hull $\wh K\ss X$ is defined in \JJ.1, and if $\wh K = K$, then $K$ is 
called {\sl Lagrangian convex} in $X$.

\Theorem{\II.8. (Contact Sets for Lag-Harmonics)}
{\sl  Let $H$ be a Lag-harmonic function on a connected open set  $X\ss \cn$. 
 Suppose $f$ is Lag-pluriharmonic (a traceless Hermitian polynomial)
with $f\leq H$ on a compact set $K\ss X$. Consider the {\rm contact set}
$$
\Sigma \ \equiv\ \{x\in K : f(x) = H(x)\}
$$
Then 
$$
\Sigma \ \ss\ \wh{\Sigma \cap \partial K.}
$$
In particular, if $K$ is Lagrangian convex in $X$, then}
$$
\Sigma \ =\ \wh{\Sigma \cap \partial K.}
$$
\pf  This follows from   [\CONT] where quite general theorems of this sort are proved.\qed

\Remark{\II.9}  Theorem \II.8  has analogues for the subequations $\cp$ and $\cp_\bbc$.
The real case was first done in [O] and [OS]. For the complex case we have that  $H$ is complex Monge-Amp\`ere
harmonic and $f$ is the real part of a complex polynomial.  Both cases are covered by [\CONT].

\Remark{\II.10}   In Section \LL\  the Dirichlet problem is solved for the  Lagrangian Monge-Amp\`ere operator 
 by the Perron method.  That is,  for a domain $\O$ and a function $\vf\in C(\bo)$, one takes the 
upper envelope of the  $u\in \cp_{\Lag}(\O)\cap \USC(\ob)$ with $u\leq \vf$ on $\bo$.  Once the theorem
is proved, one observes that smaller families will also work.  For example, one can use 
$u\in \cp_{\Lag}(\O)\cap C(\ob)$ with $u\leq \vf$ on $\bo$.

Now  in $\cn$ we have  the Lag-pluriharmonics, and one might wonder whether the solution can
be realized using   just those functions.  This is in fact the case.  There is a general
theorem which applies to all {\sl basic edge subequations},  of which $\LAG$ is one.
Details appear in [\EDGE].




\vfill\eject

\centerline{\bf  \CC. A Lagrangian Operator of Monge-Amp\`ere-Type  }
\medskip

 In this section we introduce a nonlinear partial differential operator $\mlag(D^2 u)$ 
 which defines $\partial \cp(\LAG)$ and the subequation $\cp(\LAG)$, but also has other branches.
 This puts Lagrangian potential theory in a very special category  of nonlinear equations
 whose constraint set is naturally defined  using a polynomial operator.
 
 As in the previous section, given $A\in\Sym_\bbr(\cn)$, let 
 $$
 \mu\ \equiv {\tr A\over 2},  \and \l_1, ... , \l_n\ \geq\ 0
 $$
denote the non-negative eigenvalues of the skew-hermitian part 
$
A^\sk \equiv \half(A+JAJ)$ of $A$.

\Def{\CC.1. (The Lagrangian MA-operator)}
This operator $$\mlag : \Sym_\bbr(\cn)\to \bbr$$ is defined by
$$
\mlag(A)\ \equiv\ \prod_{\pm}^{2^n} \left( \mu \pm \l_1\pm \cdots \pm \l_n\right).
\eqno{(\CC.1)}
$$

\Note {(n=1)}  In this case $\mlag(A) = \mu^2-\l^2$ where $\mu=\half \tr A$ and $\pm \l$ are the eigenvalues
of $A^\sk$. One can compute that 
$$
\mlag(A) \ =\ \det A
\eqno{(\CC.2)}
$$
is the standard Monge-Amp\`ere operator on $\bbr^2$.

\Prop{\CC.2. (n $\geq$ 2)}
{\sl $\mlag(A)$ is a polynomial of degree $2^n$ on $\Sym_\bbr(\cn)$.
In fact, after setting $\l_0\equiv \mu$, it is a symmetric polynomial of degree $2^{n-1}$
in the variables $(\l_0^2, \l_1^2, ... , \l_n^2)$.  (Note that $ \l_1^2, ... , \l_n^2$ are the eigenvalues
of $\left(A^\sk\right)^2$ where each $\l_j^2$ occurs with multiplicity 2).
}

\pf
$$
\eqalign
{
\mlag(A) \ &\equiv\ \prod_{\pm}^{2^n} \left( \mu \pm \l_1\pm \cdots \pm \l_n\right)  \cr
&=\ \prod_{\pm}^{2^{n-1}} \left[  \left( \mu \pm \l_1\pm \cdots \pm \l_{n-1}\right) + \l_n\right]
 \prod_{\pm}^{2^{n-1}} \left[  \left( \mu \pm \l_1\pm \cdots \pm \l_{n-1}\right) - \l_n\right]   \cr
 &=\  \prod_{\pm}^{2^{n-1}} \left[  \left( \mu \pm \l_1\pm \cdots \pm \l_{n-1}\right)^2 - \l_n^2\right].
}
$$
This last expression is fixed by $\mu\mapsto -\mu$ and by interchanging $\mu$ and $\l_1$.
Since the first expression is fixed by the permutation group $\pi_n$ acting on $(\l_1, ... ,\l_n)$
and by $\bbz_2^n$ acting on $(\l_1, ... , \l_n)$ by $\pm$, this proves the proposition.\qed
\medskip
 Proposition \CC.2 gives the polynomial $\mlag(A)$ a particularly
nice structure.   We illustrate this in the first two cases.

\medskip
\noindent
{\bf Example. (n=2).}
$$ 
\eqalign
{
\mlag(A)\ &=\ \l_0^4+\l_1^4+\l_2^4 -2\left( \l_0^2\l_1^2+\l_0^2\l_2^2+\l_1^2\l_2^2  \right)   \cr
&= \ \mu^2 -2(\l_1^2+\l_2^2)\mu^2 + (\l_1^2-\l_2^2)^2
}
$$

\medskip
\noindent
{\bf Example. (n=3).}
$$ 
\eqalign
{
\mlag(A)\ &=\ \left[\l_0^4+\l_1^4+\l_2^4+\l_3^4 -  2\left( \l_0^2\l_1^2+\l_0^2\l_2^2+\l_0^2\l_3^2
+\l_1^2\l_2^2+\l_1^2\l_3^2+\l_2^2\l_3^2  \right)\right]^2  \cr
&\qquad\qquad\qquad\qquad \qquad        -(8\l_0\l_1\l_2\l_3)^2      \cr
&=\ \left[\mu^4+\l_1^4+\l_2^4+\l_3^4 -  2 \mu^2\left(\l_1^2+\l_2^2+\l_3^2\right) -2
\left(\l_1^2\l_2^2+\l_1^2\l_3^2+\l_2^2\l_3^2  \right)\right]^2  \cr
&\qquad\qquad\qquad\qquad \qquad        -(8\mu\l_1\l_2\l_3)^2
}
$$

\medskip

Central to the discussion of operators of this type is G\aa rding's theory of hyperbolic polynomials.
(We refer the reader to  ([G\aa], [\HLG], [\HLGG]) for details.) To begin this discussion we need the following.

\Lemma{\CC.3}
$$
\mlag\left({t\over n} \id + A\right)\ =\  \prod_{\pm}^{2^n} \left(t+ \mu \pm \l_1\pm \cdots \pm \l_n\right).
\eqno{(\CC.3)}
$$

\pf Note that $\half\tr({t\over n} \id+A) = t+{\tr A\over 2}$,
and ${t\over n} \id+A$ has the same skew-hermitian part as $A$.\qed

\Cor{\CC.4}
{\sl
$\mlag$ is a G\aa rding polynomial hyperbolic in the direction ${ 1 \over n} \id$ 
with $\mlag({1\over n}\id)=1$.  The G\aa rding eigenvalues of $\mlag({t\over n}\id+A)$,
which by definition are the negatives of the roots of this polynomial in $t$, are}
$$
\Lambda_{\pm\pm\cdots\pm}(A) \ \equiv \ \mu\pm \l_1\pm \cdots\pm \l_n.
\eqno{(\CC.4)}
$$

\Prop{\CC.5}
{\sl
The subequation $\cp(\LAG)$ is exactly the 
 closed G\aa rding cone  defined by the inequalities $\Lambda_{\pm\pm\cdots\pm}(A)\geq0$.
}

\pf
With $\l_1, ... , \l_n \geq0$ denoting the non-negative eigenvalues of $A^\sk$, the minimum
G\aa rding eigenvalue is $\L_{\rm min}(A) = \mu-\l_1-\cdots -\l_n$. Theorem \BB.6 states that 
$\L_{\rm min}\geq0$ defines $\cp(\LAG)$.\qed

\Cor{\CC.6}
{\sl
The G\aa rding polymonial $\mlag$ is a G\aa rding/Dirichlet operator in the terminology
of [\HLG], that is, the associated closed G\aa rding cone $\cp(\LAG)$  satisfies $\cp\ss \cp(\LAG)$.
}
\pf
One has $\cp\ss \cp(\LAG)$, since $P\geq0 \ \Rightarrow \ \tr(P\bigr|_W)\geq0$ for any subspace $W$.
\qed

\medskip

Note that (\CC.2) is not all that surprising since, when $n=1$,  $\LAG = G((1,\bbr^2)$ 
and so $\cp(\LAG)=\cp$, the subequation which defines convex functions.

\vskip .3in
\centerline
{\bf
$\mlag(D^2 u)$ as a Polynomial Differential Operatior
}
\medskip

While Proposition \CC.2 gives $\mlag(A)$ a nice structure (symmetrically intertwining $\mu^2 $ and the $\l_j^2$'s,) it somewhat obscures  $\mlag(D^2 u)$)
as a differential operator.  We now give two expressions for $\mlag(D^2 u)$ in terms of the second coordinate partial derivatives of $u$.  By Proposition \CC.2, there exists 
homogeneous symmetric polynomials $s_k(\l) \equiv s_k((A^\sk)^2)$
of degree $k$ in  $\l_1^2, ... , \l_n^2$ (degree 2$k$ in $\l_1,...,\l_n$)  for $k=1,..., 2^{n-1}$,
such that:
$$
\mlag(A)\ =\   \mu^{2^n} + s_1 \mu^{2^n-2} + s_2 \mu^{2^n-4}  + \cdots + s_{2^{n-1}}
\eqno{(\CC.5)}
$$
With an abuse of notation, if we set  
$$
s_k(u) \ \equiv\ s_k \left(\left( (D^2 u)^\sk\right)^2\right),
$$
then
$$
\mlag(D^2 u)\ =\   (\D u)^{2^n} + s_1(u)  (\D u)^{2^n-2} + s_2(u)  (\D u)^{2^n-4}  + \cdots + s_{2^{n-1}}(u)
\eqno{(\CC.5)'}
$$
Rather than focus on computing what these functions $s_k(\l)$ are explicitly, we move on to our second expression for $\mlag(D^2 u)$.

Each $s_{k}(\l)$ can be written as a polynomial expression
in the basic functions
$$
\tau_{\ell}  \ \equiv\  \l_1^{2\ell}+\cdots+ \l_n^{2\ell}
$$

Now to compute the operator $\mlag(D^2 u)$ we  recall that  $\mu = \half \D u$ and that
$\l_1^2, ... , \l_n^2$ are the eigenvalues (with multiplicity 2) of $[(D^2 u)^\sk]^2$.  
Given a symmetric polynomial  $\s(\l)$ in the $\l_j^2$ we set 
$$
\s(u) \ \equiv \ \s(\l_1^2, ... , \l_n^2).
$$ 
In particular, we have
$$
\tau_{\ell}(u) \ =\  \half\tr\left\{ [(D^2 u)^\sk]^{2\ell}\right\}
$$
Calculation shows the following.
\smallskip
\noindent
{\bf Example. (n=2).} 
$$ 
\mlag(A)\ =\ \mu^4 - 2\tau_1\mu^2 +  (2\tau_2 -\tau_1^2)
$$
$$
\mlag(D^2 u) \ =\ (\D u)^4 - 2 \tau_1(u) (\D u)^2 + 2\tau_2(u) -\tau_1(u)^2
$$

\smallskip
\noindent
{\bf Example. (n=3).} 
$$
\mlag(A)\ =\ =\ \mu^8 - 4\tau_1\mu^6 + (4\tau_2 +2 \tau_1^2 )\mu^4
+\smfrac 4 3 [-16\tau_3 +18\tau_2\tau_1 -5\tau_1^3] \mu^2+ (2\tau_2-\tau_1^2)^2
$$
$$
\mlag(D^2 u) \ =\ (\D u)^8 -4\tau_1(u) (\D u)^6 + \{ 4\tau_2(u) + 2 \tau_1(u)^2\}(\D u)^4+\cdots
$$

\vfill\eject

\centerline{\bf Branches of $\mlag(D^2 u)$}
\medskip
 It is a general fact  (see [\HLG]) that any operator defined by a G\aa rding/Dirichlet polynomial
 of degree $d$ on $\Sym(\bbr^N)$ has $d$ {\bf branches}, i.e., $d$ subequations 
 defined by requiring:
$$
 \L_k(A)  \ \geq \ 0,
$$
 where $\L_1(A) \leq \L_2(A) \leq\cdots\leq \L_d(A)$ are the ordered G\aa rding eigenvalues of $A$.
The smallest such subequation, defined by $\L_1(A) \equiv \L_{\rm min}(A) \geq 0$,
 is convex and is called the {\sl G\aa rding cone}. It provides the {\sl monotonicity
property}  
$$
\{\L_1(A)\geq0\} + \{\L_k(A)\geq0\} \ \ss\ \{\L_k(A)\geq0\}.
\eqno{(\CC.6)}
$$
\Def{\CC.7} By Corollary \CC.4 the subequation $\cpL$ is the G\aa rding cone associated to
the G\aa rding/Dirichlet polynomial (\CC.1) with eigenvalues given by (\CC.4).  This gives us $d=2^n$ 
branches  
$$
\cp_k(\LAG) \ \equiv \ \{A : \L_k(A)\geq0\}, \qquad  k=1,2,...
$$ 
of the equation $\mlag$, each monotone with respect to the principle branch $\cpL=\cp_1(\LAG)$
(and so, in particular, they are subequations).  The largest  branch $\cp_d(\LAG)$ is the dual subequation
$\wt{\cp(\LAG)}$.

\medskip

This monotonicity of the branches with respect to $\cpL$, has interesting consequences.
For example, it implies the following removable singularity result.

\Theorem {\CC.8.\ ([\RSing])} 
{\sl
Let $\O\ss\rn$ be a domain
and $\Sigma\ss \O$ a closed subset of locally finite Hausdorff $(n-2)$-measure.
 Then any function $u$ which is  $\cp_k(\LAG)$-subharmonic
on $\O-\Sigma$ and is locally bounded above at points of $\Sigma$, extends to a $\cp_k(\LAG)$-subharmonic function on all of $\O$.  Similarly, if $u$ is $\cp_k(\LAG)$-harmonic on $\O-\Sigma$ and continuous
on $\O$, then $u$ is $\cp_k(\LAG)$-harmonic on $\O$.
}

\Note{\CC.9.  (Unitary Invariance)}  The operator $\mlag$ on $\Sym_\bbr(\bbc^n)$ is U$(n)$-invariant, i.e.,
$$
\mlag(gA g^{-1})=\mlag(A)\quad {\rm for\ all\ }\ g\in {\rm U}(n) \ \ {\rm and } \ \ A\in \Sym_\bbr(\bbc^n).
$$
This follows rather straightforwardly form the definition of $\mlag$.  It follows that each of the subequations
$$
\cp_k(\LAG)\quad {\rm is\ \ } {\rm U}(n) \ {\rm invariant, \ for\ \ } 1\leq k\leq 2^n.
$$
This fact is important because it shows that the subequations $\cp_k(\LAG)$ make sense on
any Gromov manifold (see \S \FF).


\vfill\eject

\centerline{\bf  \DD. Two Intrinsic Definitions of the Lagrangian (Monge-Amp\`ere) Operator}
\bigskip

In this section we present two constructions of the operator $\mlag(D^2 u)$.
Each gives us an intrinsic construction of the operator when we pass from $\cn$ to general
Gromov manifolds.

\vskip .3in
\centerline{\bf  A.\  Constructing   $\mlag$ as the Determinant of a Derivation.}
\medskip

Consider $\bbc^n = (\bbr^{2n}, J)$ as above.  Then any  symmetric  matrix $H\in \Sym(\bbr^{2n})$
prolongs to  a map 
$$
D_H : \L^n \bbr^{2n} \arr  \L^n \bbr^{2n}
$$
 as a derivation.

\medskip
\noindent
{\bf Observation \DD.1.}
{\sl
Suppose that $H$ can be diagonalized by a hermitian orthonormal basis
$e_1, Je_1, ... , e_n, Je_n$.  Then $D_H$ preserves the subspace $S_H\ss \L^n \bbr^{2n}$  spanned by
the $2^n$ vectors
$$
\xi_{\pm\pm\cdots\pm} \ \equiv\ (e_1 \ {\rm or} \ Je_1) \wedge (e_2 \ {\rm or} \ Je_2) \wee (e_n \ {\rm or} \ Je_n)
$$
In fact these vectors are eigenvectors of $D_H$.
}

\medskip

We can now apply this observation to the Lagrangian part $H \equiv A^\Lag = {1\over 2n}(\tr A)\id+A^\sk$
of any matrix $A\in \Sym_\bbr(\cn)$, and straightforward calculation shows that
$$
\mlag(A) \ =\ \det\left\{ D_{A^{\Lag}} \bigr|_{S_{A^\Lag}} \right\}.
\eqno{(\DD.1)}
$$
With notation as in Definition \DD.1, the eigenvector $\x_{\pm\pm\cdots\pm}$ of $D_H$ has eigenvalue
precisely $\mu\pm\l_1 \pm\l_2 \pm\cdots\pm\l_n$.)
In particular,
$$
\mlag(A) \ \ {\rm is\ a\ factor\ of\ the\ determinant\ of\ \  }   D_{A^{\Lag}}  
\ \ {\rm acting \ on \ \ } \L^n{\bbr^{2n}}.
\eqno{(\DD.2)}
$$
We note that  although (\DD.1) does not provide an intrinsic definition of $\mlag(A)$  (since
$S_{A^\Lag}$ depends on $A$), the factor $\mlag(A)$ in (\DD.2) is intrinsic.

The assertion (\DD.2) can be somewhat generalized by considering the subspace
$$
\L^n_{\rm prim} 
\ = \ \{\vf \in \L^n : \o \wedge \vf=0\}
$$

\Lemma{\DD.2} 
{\sl
The subspace $\L^n_{\rm prim}  \ss \L^n\bbr^{2n}$ is invariant under $D_{A^{\Lag}}$
for any $A\in  \Sym_\bbr(\cn)$.  Thus,  in particular,}
$$
\mlag(A) \ \ {\rm is\ a\ factor\ of\ the\ determinant\ of\ \  }   D_{A^{\Lag}}  
\ \ {\rm acting \ on \ \ } \L^n_{\rm prim} .
\eqno{(\DD.3)}
$$

\pf 
Let $e_1, Je_1,..., e_n, Je_n$ be a hermitian orthonormal basis which diagonalizes
$A^\sk$.  The since $\o= e_1\wedge Je_1 +\cdots + e_n\wedge Je_n$, we have
$D_{A^\sk} (\o)=0$, and so $D_{A^{\Lag}} (\o)=  {\tr(A)\over n} \o$.
We can assume $\tr(A)\neq0$ by adding a multiple of the identity if necessary.
Then if  $\psi\in \L^n_{\rm prim}$, we have
$$
0\ =\ 
D_{A^{\Lag}} (\o \wedge \psi) 
\ =\ D_{A^{\Lag}}(\o) \wedge \psi + \o\wedge D_{A^{\Lag}}(\psi)
\ =\  \o\wedge D_{A^{\Lag}}(\psi),
$$
and so $ D_{A^{\Lag}}(\psi)  \in \L^n_{\rm prim}$ .\qed

\vfill\eject
\centerline{\bf  B.\  Constructing   $\mlag$ as the Determinant of a Spinor Endomorphism.}
\medskip

In this subsection we present a construction of the operator $\mlag(D^2 u)$ in terms of spinors 
-- more specifically in terms of the irreducible representations of  Spin$^c_{2n}$.  This gives a 
second intrinsic construction of the operator, which is useful when we pass from $\cn$ to general
Gromov manifolds.  We break the construction into two steps.

The first step is to define a natural algebraic map
$$
\Phi : \Sym_\bbr(\cn) \ \arr\ \L^2_\bbr \cn.
\eqno{(\DD.4)}
$$
Given $A\in \Sym_\bbr(\cn)$, abbreviate 
$$
B \ \equiv\ A^\sk \ =\ \half(A+JAJ).
$$
Since $Be=\l e$ implies $BJe=-\l Je$, the square $B^2\geq 0$ has a unique positive
square root $E\equiv \sqrt{B^2}$ which satisfies $Ee=|\l| e$ and $EJe = |\l | Je$.  Note that
$JE=EJ$ and that $E$ is hermitian symmetric. Hence,
$$
EJ \in {\rm Skew}_\bbr(\cn)\ \equiv\ \L^2_\bbr \cn, \qquad {\rm i.e.,}\ \ (EJ)^t = - EJ.
\eqno{(\DD.5)}
$$

\Def{\DD.3}  Adopting these notations the map $\Phi$ is defined by:
\medskip

(a)\ \ $\Phi(A) \ \equiv\  EJ\ =\ \sqrt{B^2} J, \ \ B\ \equiv\ \half(A+JAJ)$, 
\medskip

\noindent  or using Lemma \BB.1,
\medskip
 
(b)\ \ $\Phi(A) \ \equiv\  \l_1 e_1\wedge Je_1 + \cdots + \l_n e_n\wedge Je_n$.
\medskip

To see that (a) and (b) are equal, apply the canonical form for 
$B\in \Herm^\sym(\cn)$ (with $\l_1\geq \cdots\geq \l_n\geq 0)$, and  note that
$$
E\ \equiv \ \l_1(P_{e_1}+P_{Je_1}) + \cdots + \l_n(P_{e_n}+P_{Je_n}),
\eqno{(\DD.6)}
$$
and second that
$$
EJ \ =\ \l_1(Je_1 \circ e_1 - e_1 \circ Je_1) + \cdots + \l_n(Je_n \circ e_n - e_n \circ Je_n)
\eqno{(\DD.7)}
$$

\def\bbB{{\bf B}}\def\bcl{{\bf Cl}_{2n}}

\Remark {\DD.4}
Consider the blocking $\cn \equiv W\oplus JW$ defined by (\BB.9).
First $J\equiv \left( \matrix{0 & -I \cr I & 0}\right)$, so that 
$$
B\equiv \left( \matrix{\l & 0 \cr 0 & - \l}\right), \qquad
E \equiv \left( \matrix{\l & 0 \cr 0 & \l}\right), \quad {\rm and}\qquad
EJ \equiv  \left( \matrix{0 & -\l \cr \l & 0}\right).
\eqno{(\DD.8)}
$$
Note that $\l \geq 0$ is diagonal.  Also note that $B = B^+ - B^-$ with
$$
B^+ \ \equiv \ \half(B+E) \ =\ 
\left( \matrix {\l &0 \cr 0 & 0 }\right) \ \geq\ 0
\and 
B^- \ \equiv \ \half(B-E) \ =\ 
\left( \matrix {0 &0 \cr 0 & \l }\right) \ \geq\ 0
$$
the positive and negative parts of $B$.

For the second step adopt the notation $A\in \Sym_\bbr(\cn)$ and $B\equiv A^\sk$ as before.
Now using the first step,  set $\wt \bbB \equiv \Phi(A)$ and 
 consider  $\wt \bbB$ as an element in the Clifford algebra
$$
\wt \bbB \ \in \ \L^2 \bbr^{2n} \ \ss\ {\rm Cl}_{2n}.
$$
Recall that
$$
\wt \bbB \ =\ \sum_{k=1}^n \l_k e_k Je_k.
$$

Consider now the (unique) irreducible complex representation $\SS \cong \bbc^{2^n}$ of 
Cl$_{2n}$. This extends naturally to a representation
of $\bcl = {\rm Cl}_{2n}\otimes_{\bbr} \bbc$.  Note that the elements
$$
i e_k Je_k \in \bcl  \quad{\rm satisfy}\quad (i e_k Je_k)^2\ =\ 1.
$$
We set
$$
\pi^+_k \equiv \half(1+i e_k Je_k)
\and
\pi^-_k \equiv \half(1-i e_k Je_k).
$$
Then
$$
(\pi^+_k)^2\ =\ \pi^+_k, \quad (\pi^-_k)^2\ =\ \pi^-_k
\and \pi^+_k+\pi^-_k \ =\ 1.
$$
Also we have that the $\pi^\pm_k $ commute with all the $ \pi^\pm_\ell$ for $k\neq \ell$ and 
$$
\pi^\pm_k e_k \ =\ -e_k \pi^\pm_k\qquad{\rm for\ all \ \ } k.
$$
Thus Clifford multiplication by the projectors $\pi^\pm_k$ decompose the Spinor space
$$
\SS\ = \ \equiv\ \pi^+_k \SS \ \oplus \ \pi^-_k \SS
$$
into the $+1$ and $-1$ eigenspaces under multiplication by $ie_k Je_k$. These two spaces are isomorphic
under multiplication by $e_k$ and they are invariant under multiplication by $ie_\ell J e_\ell$ for all
$\ell\neq k$.  Hence we get a decomposition
$$
\SS\ =\ \bigoplus_{{\rm all}\ \pm} \SS_{\pm \,\pm\, \cdots \pm}
\qquad{\rm where}\qquad  \SS_{\pm \,\pm\, \cdots \pm} \ \equiv \ \pi^\pm_1 \pi^\pm_2 \cdots \pi^\pm_n \SS
\ =\ \bigcap_{k=1}^n \pi_k^{\pm} \SS.
$$
Now each $\SS_{\pm \,\pm\, \cdots \pm}$ has dimension one  and is 
an eigenspace for each $i e_k Je_k$. This follows because the $\pi^\pm_k$ all commute and 
$\pi_k^+ e_k = - e_k \pi_k^-$,  and so one can do a dimension count (see [LM, page 43 ff.]).
 In particular each 
$$
\SS_{\pm \,\pm\, \cdots \pm} \ \ {\rm
is \  an \  eigenspace\  for \ \ } i \wt \bbB\ \ {\rm with\ eigenvalue\ \  }
\pm \l_1\pm \l_2 \pm\  \cdots\ \pm\l_n.
\eqno{(\DD.9)}
$$
Thus we have:
\Prop {\DD.5}
{\sl As an endomorphism of $\SS$ by Clifford multiplication, the element $i \wt \bbB$ has
$$
\det(i \wt \bbB) \ =\ \prod_{{\rm all}\ \pm} \left( \pm \l_1\pm\l_2\pm \cdots\pm\l_n\right)
$$
Furthermore, for any real number $\mu$, under Clifford multiplication, the element $ \mu 1 +i \wt \bbB$ has
$$
\det(\mu 1+i \wt \bbB) \ =\ \prod_{{\rm all}\ \pm} \left(\mu \pm \l_1\pm\l_2\pm \cdots\pm\l_n\right)
$$
}

\vfill\eject


\vfill\eject

\centerline {\bf \FF. Lagrangian Potential Theory on Gromov Manifolds.}   
\medskip

In this section  we carry the previous discussion over to a  general context which is
relevant to symplectic geometry.

\Def{\FF.1}  By a  {\bf Gromov manifold} we mean a triple  $(X,\o,J)$ 
where  $(X,\o)$ is a symplectic manifold and $J$ is an almost complex structure  on $X$ satisfying the
 conditions:
$$
\o(v,w) \ =\  \o(Jv, Jw) \and \o(Jv, v)\ >\ 0
\eqno{(\FF.1)}
$$
for all non-zero tangent vectors $v$ and $w$ at each point of $X$.
On such a   manifold there is a natural riemannian metric $\bra {\cdot}{\cdot}$ defined
by
$$
\bra {v}{w} \ \equiv \ \o(Jv,w) \qquad{\rm with}\quad \bra{Jv}{Jw} \ =\ \bra vw.
\eqno{(\FF.2)}
$$

\Remark {\FF.2}  It is a result of Gromov (see [G], [MS, p??]) that any compact
symplectic manifold admits an almost complex structure $J$ satisfying (\FF.1).

\Remark {\FF.3}
Note that on a Gromov manifold there is a well defined notion of Lagrangian submanifold
defined as in (\BB.1).
There are also local $J$-holomorphic curves passing through any point with any
prescribed complex tangent [NW]. One easily checks that if $F:X\to X$ is a symplectomorphism,
(a diffeomorphism with $F^*\o=\o$),  then:

\medskip
(a) \ \ $F$ maps Lagrangian submanifolds to Lagrangian submanifolds,

\medskip
(b) \ \ The transported almost complex structure $\wt J \equiv F_* \circ J\circ (F^{-1})_*$
again satisfies

\qquad  \ \ the conditions (\FF.1),

\medskip
(c) \ \ $F$ maps $J$-holomorphic curves to $\wt J$-holomorphic curves.

\medskip

We now recall that using the riemannian metric, one can define a  Hessian (or second derivative)
which will allow us to carry over the foregoing material to Gromov manifolds.

\Def{\FF.4}   Given a smooth function $u\in C^\infty(X)$, the
{\bf hessian}  of $u$ is a section of $\Sym(T^*X)$ defined on vector fields
$v,w$ by
$$
(\Hess f)(v,w)\ \equiv\ vwf-  (\nabla_vw)f
\eqno{(\FF.3)}
$$
where $\nabla$ denotes the Levi-Civita connection for the metric $\bra {\cdot}{\cdot}$.
\medskip

We note that this Hessian gives a canonical splitting of the 2-jet bundle of $X$:
$$
J^2(X) \ =\ \bbr\oplus T^*X \oplus \Sym(T^*X),
\eqno{(\FF.4)}
$$
(see [\DDR] or [\SURVEY]).
 Furthermore, using the metric and $J$ we identify $$T^*X\ \cong\  TX$$ 
 and obtain  the decomposition  
$$
\Sym(T^*X) \ =\ \bbr \oplus \Herm^\sym_0(TX) \oplus \Herm^\sk(TX)
\eqno{(\FF.5)}
$$
corresponding exactly  to the decomposition (3.3).

Observe now that there is a well defined subbundle $\LAG\ss G(n, TX)$ of the
Grassmann bundle of tangent $n$-planes, which consists of the Lagrangian 
tangent $n$-planes.  This bundle is invariant under the group of symplectomorphisms
of $(X,\o)$.
It embeds naturally into the bundle $\Sym(TX)$ by associating to $W$, the orthogonal 
projection of $T_xX$ onto $W$.

With this said it should be clear that all of the algebraic considerations of the 
previous sections 2 and 3  apply in the obvious way to this context.  In particular,
we have a well-defined subequation of the 2-jet bundle:
$$
\cp(\LAG) \ \equiv\  \left \{ J^2(u) :\tr\left( \Hess(u) \bigr|_W\right) \geq 0\ \forall\ W\in\LAG\right\}
\ \ss\ J^2(X),
$$
as well as its dual subequation defined fibrewise by $\cpt(\LAG) \equiv - (\sim \Int \cp(\LAG))$,
and the equation $\partial \cp(\LAG) = \cp(\LAG)\cap (-\cpt(\LAG))$.

We can also carry over the definitions from Section \BB. The primary two are:

\Def{\FF.5} A   function $u\in C^\infty(X)$ is {\bf Lagrangian plurisubharmonic} if $J^2_x(u)
\in \cp(\LAG)$ for all $x\in X$, or equivalently, if 
$$
\tr\left( \Hess(u) \bigr|_W\right) \geq 0\ \forall\ W\in \LAG.
$$
If, in addition, $-J^2_x(u) \in \cpt(\LAG)$ for all $x$, then $u$ is called {\bf   Lagrangian harmonic}.
Equivalently, $u$ is Lagrangian harmonic if and only if  $J_x^2(u) \in \partial \cp(\LAG)$ for all $x$.
\medskip

All of the  definitions in Section 3 extend to upper semi-continuous functions. The notion of a viscosity test
function (Definition \BB.5) carries over directly to manifolds, and we have the following.

\Def{\FF.6}  A function $u\in\USC(X)$ is {\bf  Lagrangian plurisubharmonic} if for 
each $x\in X$ and each test function $\vf$ for $u$ at $x$,  one has $J^2_x(\vf) \in \cp(\LAG)$.
If, in addition, for 
each $x\in X$ and each test function $\psi$ for $-u$ at $x$,  one has $J^2_x(\psi) \in \cpt(\LAG)$,
then $u$ is called {\bf Lagrangian harmonic}.
(This additional condition by itself defines the notion of {\bf dually Lagrangian plurisubharmonic} for 
$v\equiv -u$.)

\Theorem{\FF.7}  {\sl  If $M$ is a Lagrangian submanifold which is also a minimal submanifold,
then the restriction $u\bigr|_M$ of a Lagrangian plurisubharmonic  function $u$ to $M$ is subharmonic.
}

\pf  If $u$ is $C^\infty$ it follows from the equation (see Proposition 2.10 in [\PTCG])
$$
\D\left(u\bigr|_M\right)\ =\ \tr_{TM} \{ \Hess u \} - H_M(u)
$$
where $H_M$ is the mean curvature vector of $M$ and $\D$ is the intrinsic
Laplace-Beltrami operator on $M$ with respect to the induced metric.
For $u\in \USC(X)$ one must use the Restriction Theorem [\REST, Thm. 6.4].
\qed
\medskip

Note that (\BBB.12a) easily carries over.

\Lemma{\FF.8}  {\sl  If $u\in \USC(X)$ is Lagrangian plurisubharmonic, then $u$ is subharmonic.}
\pf Let $\vf$ be a test function for $u$ at $x\in X$, and choose a Lagrangian plane $W\ss T_xX$.
Then $\tr\{\Hess_x\vf\bigr|_W\}\geq0$ and  $\tr\{\Hess_x\vf\bigr|_{W^\perp}\}\geq0$
since $W^\perp$ is also Lagrangian.  Hence, $\D(\vf)_x = \tr\{\Hess_x\vf\bigr|_W\}
+ \tr\{\Hess_x\vf\bigr|_{W^\perp}\} \geq0$, and we conclude that $u$ is subharmonic
in the viscosity sense. However, this notion coincides with all other notions of subharmonicity
on a riemannian manifold.
\qed

\medskip
\noindent
{\bf Note \FF.9.}  In defining  our Lagrangian plurisubharmonic functions here, we have used the riemannian
Hessian.  That same Hessian could also be used to define complex plurisubharmonic functions on
our almost complex manifold.  However, these complex psh functions 
do not always agree with the usual intrinsic ones,
namely those whose restrictions to  (pseudo) holomorphic curves are subharmonic (see Example 9.5
in [\ACM]).  However, when the K\"ahler form of the almost complex structure is $d$-closed,
as it is here,  the two notions of complex psh functions do agree [\ACM, Thm. 9.1].

\vskip .3in


\centerline {\bf \HH. The Lagrangian Monge-Amp\`ere Operator  on Gromov Manifolds.}   
\medskip

On any Gromov manifold $(X, J, \o)$ (in fact on any almost complex hermitian manifold) there is 
a well-defined Lagrangian Monge-Amp\`ere Operator $\mlag(\Hess\, u)$ defined exactly as in the
euclidean case, but with $D^2 u$ replaced by the riemannian Hessian $\Hess\,u$.
As before this operator has $2^n$ branches where $2n = \dim_\bbr(X)$. This operator can be defined intrinsically
in two different ways.

For the first we consider the derivation 
$$
D \ \equiv\ D_{\Hess (u)^\Lag}
$$
acting of the bundle $\L^n T^*X$.   Then as in (\DD.2) we see that 
$$
\mlag(\Hess\, u) \   {\rm is\  a\   factor\ of\ }\ \det \left \{ D \right\}.
$$
In fact as noted there we can restrict $D$ to the 
bundle of primitive $n$-forms 
$$
\L^n_{\rm prim}(X) \ \ss\ \L^n T^*X
$$
defined as the kernel of exterior multiplication $\o\wedge : \L^n T^*X \to \L^{n+2} T^*X$
by the 2-form $\o$.  Then,  from Lemma \DD.2 we have
$$
\mlag(\Hess\, u) \   {\rm is\  a\   factor\ of\ }\ \det \left \{ D \bigr|_{\L^n_{\rm prim}(X) }\right\}.
$$

For the  second construction let $\SS\arr X$ be any bundle 
of complex modules over the Clifford bundle ${\bf Cl}(X) \equiv Cl(X)\otimes_\bbr \bbc$.
Assume further that $\SS$ is pointwise irreducible, i.e., dim$_\bbc(\SS) =2^n$.
Then given any function $\mu\in C^\infty(X)$ and any section $B$ of  $\Sym(T^*X)$ which is skew hermitian ($BJ\equiv -JB$ on $X$), we consider the section of the Clifford  bundle
$$
\mu 1 + i \wt \bbB \ \in\ \Gamma({\bf Cl}(X)).
$$
in the notation of Section  \DD.  Clifford multiplication gives a bundle map
$$
(\mu 1 + i \wt \bbB) : \SS \ \arr\ \SS
$$
and we can take its determinant.

Now there always exists such a bundle of irreducible complex modules for ${\bf Cl}(X)$, namely
$$
\SS \ \equiv \ \bigoplus_{q=0}^n \L^{0,q}(X) \ \cong \ \L_\bbc^* T_\bbc^*(X).
$$
The Clifford multiplication is generated by letting a real tangent vector $v$ act by 
$v\wedge - i_v $  where $i_v$ is contraction.  This gives the following.

\Theorem {\HH.1} {\sl  Let $u$ be a $C^2$-function on a Gromov manifold $X$.  Consider the section
$$
\half \D u + i {\wt {\bf H}}_{\rm skew}(u)
$$
acting by Clifford multiplication on the  bundle $\SS \ \equiv \ \bigoplus_{q=0}^n \L^{0,q}(X) $ of  $(0,q)$ forms.
Then the Lagrangian Monge-Amp\`ere operator on $X$ is the determinant of this bundle map:
$$
\mlag(\Hess\, u) \ =\ \det\left\{\half \D u + i {\wt {\bf H}}_{\rm skew}(u)\right\}
$$
}
\medskip

\noindent
{\bf Note \HH.2.}  There are other choices, we could twist $\SS$ with any complex line bundle.  This
does not change the determinant.
\medskip
\noindent
{\bf Note \HH.3.}  There is a "quasi" form of this equation.  We replace $\wt \bbB= {\wt {\bf H}}_{\rm skew}(u)$ with
$\o +\wt \bbB$.

\Remark{\HH.4. (Branches)}
Each of the branches $\cp_k(\LAG)$, $1\leq k\leq 2^n$ of the Lagrangian operator  (see \S \CC\ and Def. \CC.7)  carries over to any Gromov manifold.  This follows from the unitary invariance of  $\cp_k(\LAG)$, which follows 
from the unitary invariance of $\mlag$, as discussed in Note \CC.9.


\vskip .3in
\centerline{\bf  \JJ.  Lagrangian Pseudoconvexity.}
\medskip
 
In this section we investigate the \lag\ analogue of the concepts of  
a pseudoconvex domain  and  a total real submanifold in complex analysis.
Suppose $(X,\o,J,\langle\cdot,\cdot\rangle)$ is a non-compact, connected Gromov manifold, and denote by
$\PSHl^\infty(X)$ the cone of smooth \lagpsh\  functions on $X$.  
Recall also the succinct notation Lag-psh, etc. in Definition \BBB.5.

\Def{\JJ.1}  By the {\bf \lag\ hull} of a compact subset $K\subset X$ we mean the set
$$
\PH K \ \equiv \ \{x\in X: f(x)\leq \sup_K f \ \ {\rm for\ all\ \ } f\in \PSHl^\infty(X)\}
$$
If $\PH K = K$, then $K$ is called  {\bf (\lag) convex}.

\Theorem {\JJ.2}  {\sl The following are equivalent.
\smallskip

1) \ \ If $K\subset\subset X$, then $\PH K\subset\subset X$.
\medskip

2) \ \ There exists a smooth Lag-psh proper exhaustion function $f$ on $X$.
\medskip

3) \ \  There exists a neighborhood $N$ of $\infty$ in $X$ and a smooth function $v$ on $N$,

  \qquad   which is 
Lag-psh, such that $\lim_{x\to\infty} v(x) = \infty$.
}

\Def{\JJ.3}  A Gromov manifold $X$ satisfying the equivalent conditions of Theorem \JJ.2
is called {\bf \lag\ convex}.

\Theorem {\JJ.4}  {\sl The following are equivalent.
\smallskip

1) \ \ $K\subset\subset X\ \ \Rightarrow\ \ \PH K\subset\subset X$, and there exists $f\in C^\infty(X)$ which is strictly  Lag-psh.

\medskip

2) \ \ There exists a smooth strictly  Lag-psh proper exhaustion function   on $X$.
\medskip
}

\Def{\JJ.5}  When $X$ satisfies the equivalent conditions of Theorem \JJ.4, it
is called {\bf strictly \lag\ convex}.

\Theorem {\JJ.6}  {\sl The following are equivalent.
\smallskip

1) \ \ $K\subset\subset X\ \ \Rightarrow\ \ \PH K\subset\subset X$, and there exists $f\in C^\infty(X)$ which is strictly
 Lag-psh

\qquad  outside a compact subset of $X$.

\medskip

2) \ \ There exists a smooth   Lag-psh proper exhaustion function  on $X$  which is strict

\qquad  outside a  compact subset.\medskip

3) \ \  There exists a neighborhood $N$ of $\infty$ in $X$ and a smooth function $v$ on $N$,

  \qquad   which is strictly
Lag-psh, such that $\lim_{x\to\infty} v(x) = \infty$.

 }

\Def{\JJ.7} When $X$ satisfies the equivalent conditions of Theorem \JJ.6, it
is called {\bf strictly \lag\ convex at infinity}.

\medskip
 Theorems \JJ.2, \JJ.4 and \JJ.6 are proved  (in greater generality) 
in \S 4 of [\GPandC].

\vskip.3in
\centerline{\bf Cores}
\medskip

Given a function $f\in\PSHl(X)$, consider the open set
$$
S(f) \ \equiv\ \{x\in X : f \ \ {\rm is\  strictly\  \lag\ \psh\ at\ } \ x\}
$$
and the closed set
$$
W(f)\ \equiv \ X-S(f).
$$
Note that
$$
W(\l f+\mu g) \ =\ W(f) \cap W(g)
$$
for $f,g\in \PSHl(X)$ and $\l, \mu>0$.

\Def{\JJ.8}  The {\bf core} of $X$ is defined to be the intersection
$$
\Core(X)\ \equiv\ \bigcap W(f)
$$
over all $f\in \PSHl(X)$.  The {\bf inner core} is defined to be the set 
InnerCore$(X)$ of points $x$ for which there exists $y\ne x$ with the property
that $f(y) = f(x)$ for all $f\in \PSHl(X)$.  \medskip

Arguing exactly as in [\PTCG] shows the following:
$$
{\rm InnerCore}(X)\subseteq\Core(X)
\eqno{(\JJ.1)}
$$
$$
{\rm Every\ compact\  minimal\ Lagrangian\ submanifold\ is\ contained\ in\ } \Core(X)
\eqno{(\JJ.2)}
$$

\Theorem{\JJ.9} {\sl  Suppose $X$ is \lag\ convex. Then $\Core(X)$ is compact if and only if 
$X$ is strictly \lag\ convex at infinity.  Furthermore, $\Core(X)=\emptyset$
if and only if  $X$ is strictly  \lag\ convex.}

\vskip .3in
\centerline{\bf Free Submanifolds}
\medskip

In analogue with the concept of a totally real submanifold (one free of any complex tangent lines) in complex analysis,
we introduce the following.

\Def{\JJ.10}  A submanifold $M\subset X$ is said to be {\bf (\lag) free} if its tangent spaces
contain no \lag\ $n$-planes.

\Ex {\JJ.11}  Note that any submanifold of dimension $<n$ is automatically free. Note also that any
(almost) complex submanifold is also free.

As in [\PTCG]   free submanifolds can be used to construct huge families of strictly
\lag\ convex spaces.  We begin with the following observation.

\Theorem{\JJ.12}  {\sl Suppose $X$ is strictly \lag\ convex and of dimension 2$n$.  Then $X$
has the homotopy-type of a CW-complex of dimension $\leq 2n-2$.}

\pf Let $f:X\to \bbr^+$ be a strictly \lagpsh\ proper exhaustion function. By perturbing we may assume that  $f$ has non-degenerate critical points. The theorem follows if we show that each critical point
has Morse index $\leq 2n-2$ (cf. [M]). If this fails, then there is a critical point $x$ at which
 $\Hess_x f$ has at least  $2n-1$ negative eigenvalues.  This means there exists a subspace
$W\subset T_xX$  of dimension $\geq 2n-1$ with  $\Hess_x f\bigr|_W <0$.   However, any such $W$  contains a \lag\  $n$-plane $W$, and since $f$ is strictly \lagpsh, we must have $\tr_W \Hess_x f >0$, a contradiction.
\qed
\medskip

The following results are proved by an adaptation of [HW$_{1,2}$] 
and Theorems 6.4 and 6.5 in  [\PTCG].

\Theorem {\JJ.13}  {\sl Suppose $M$ is a closed submanifold of $X$ and let 
$f_M(x) \equiv {1\over2}{\rm dist}(x,M)^2$ denote half the square of the distance to $M$.  Then
$M$ is \lag\ free if and only if 
the function $f_M$ is strictly  Lag-psh at each point in $M$ (and hence in a neighborhood of $M$).
}

\Theorem{\JJ.14}  {\sl  Suppose $M$ is any  \lag\ free submanifold of $M$.
  Then there exists a fundamental neighborhood system $\cf(M)$ of $M$ such that:
  
  \smallskip
  
  (a)\  $M$ is a deformation retract of each $U\in \cf(M)$.

  \smallskip
  
  (b)\  Each neighborhood $U\in \cf(M)$ is strictly \lag\ convex.

  \smallskip
  
  (c)\  \ ${\PSHl}(V)$ is dense in  ${\PSHl}(U)$ if $U\subset V$ and $V,U\in\cf(M)$.

  \smallskip
  
  (d)\  Each  compact set $K\subset M$ is ${\PSHl}(U)$-convex for each $U\in \cf(M)$.
  }


\vskip .3in

\centerline{\bf  \KK.  Lagrangian Boundary Convexity.}
\medskip

Suppose $\O\subset\subset X$ is an open set with smooth boundary in a non-compact Gromov manifold $X$.    
The three global conditions in Definitions \JJ.3, \JJ.5 and \JJ.7 can be applied to the domain $\O$ since it is
also a Gromov manifold.  In this section we introduce local conditions on  its boundary $\bo$ which are of
a companion nature, and we prove a local to global result (Theorem \KK.5) in the vein of the Levi problem in 
complex analysis.

A \lag\ $n$-plane $W$ at a point  $x\in \bo$ will be called a  {\bf tangential Lagrangian} plane if 
$W\subset T_x\bo$.  Let $LAG(\bo)$ denote the set of all such planes at all points of $\bo$.

  \Def{\KK.1}  Suppose that $\rho$ is a {\bf defining function} for  $\bo$,
  that is, $\rho$ is a smooth function defined on a neighborhood 
  of $\overline{\Omega}$ with $\Omega =\{x: \rho(x)<0\}$ and  $\nabla \rho \neq 0$ on $\partial \Omega$.
  If
  $$
 \tr_W \Hess \rho \ \geq\ 0 \ \ \  {\rm for\ all\  } W \in LAG(\bo),
   \eqno{(\KK.1)} 
 $$
  then  $\partial \Omega$ is  called {\bf  \lag\ convex}.
  If the inequality in (\KK.1) is strict for all $W \in LAG(\bo)$, then 
  $\partial \Omega$ is  called {\bf strictly \lag\ convex}.  If 
  $  \tr_W \Hess \rho  =0$ for all $W$ as in (\KK.1), then    $\partial \Omega$ is {\bf \lag\ flat}.
  
  Each of these conditions is a local condition on $\bo$, independent of the choice of $\rho$.

  \Lemma {\KK.2} {\sl  Each of the three conditions in Definition \KK.1  is independent of the choice of defining function $\rho$.
  In fact, if $\overline{\rho}=f\rho$ is another choice with $f>0$ on $\partial \Omega$, then on $\bo$ }
  $$
 \tr_W \Hess \overline{\rho}\ =\ f \cdot \tr_W \Hess \rho  \ \ \ {\rm for \ all \ }\ \ W\in LAG(\bo)
    \eqno{(\KK.2)}
    $$
    
  \pf
  Note that  from (\FF.3) we have 
  $$\eqalign{
  \{\Hess (f\rho)\}(v,w) \  &= \ vw(f\rho)-(\nabla_vw)(f\rho)  \cr
  &= \  f\Hess ( \rho)(v,w)  +(vf)(w\rho)+(wf)(v\rho)
 +\rho \Hess ( f)(v,w) 
  }
  $$
  Since $\rho=0$ on $\bo$ and $v\rho=0$ for all  $v\in T(\bo)$, the assertion follows. 
  \qed

  \Lemma {\KK.3} {\sl Suppose $\rho$ is a smooth real-valued function 
  on $X$, and $\psi:\bbr\to\bbr$ is smooth on the image of $\rho$.
  Then
    $$
\tr_W\Hess \psi(\rho)  \ =\ \psi'(\rho)\tr_W \Hess\rho +   \psi''(\rho)|\nabla\rho\hk W|^2 
  \eqno{(\KK.3)}
  $$
 for all oriented tangent $p$-planes $W$.}
 
 \pf
 We first calculate that
 $
 \Hess \psi(\rho)  = \psi'(\rho)  \Hess\rho +   \psi''(\rho)\nabla \rho \circ \nabla \rho
 $
and then note that $\tr_W( \nabla \rho \circ \nabla \rho )= |\nabla\rho\hk v|^2 $.\qed

\Cor{\KK.4}  {\sl  With $\delta = -\rho$ and $\rho <0$, one has }
$$
\tr_W \Hess(-\log \delta)\ =\ {1\over \delta}\tr_W \Hess \rho + {1\over \delta^2}|\nabla \rho\hk W|^2
\eqno{(\KK.4)}
$$
\pf 
Take $\psi(t)=-\log(-t)$ for $t<0$, and note that $\psi'(t)=-1/t $ and  $\psi''(t)=1/t ^2$, so that $\psi'(\rho)=
 1/\delta$  and $\psi''(\rho)= 1/\delta^2$.\qed
\medskip

   We now come to the main result of this section, going from local to global.

  \Theorem{\KK.5}   {\sl   Let $\Omega\subset\subset X$ be a compact domain with   smooth,
    strictly \lag\ convex boundary.          Suppose $\rho$   is an arbitrary  defining function for $\bo$ 
    and let  $\delta \equiv -\rho$ be the corresponding interior  ``distance function''  to $\bo$.
  Then $-\log \,\delta$ is strictly Lag-psh    outside a compact     subset of $\O$. Thus, in particular, the domain $\Omega$ is strictly \lag\ convex at infinity.   }
  $ \over $
  
 \pf   
At each point $x\in\O$
 near $\bo$, we have that equation (\KK.4) holds for all $n$-planes $W$.
   Note that at $x\in \bo$,  $ |\nabla\rho\hk W|^2$ vanishes if and only if $W$ is tangential to $\bo$.    For notational convenience we set 
 $$
 \cos^2\theta(W) \ =\ {|\nabla\rho\hk W|^2\over |\nabla \rho|^2} \ =\ \langle P_{\span \nabla \rho}, P_{\span W} \rangle
$$
Then
the inequality $|\cos\theta|<\epsilon$ defines a fundamental neighborhood system for $G(p,T\bo)\subset G(p,TX)$.  By restriction     $|\cos\theta|<\epsilon$ defines a fundamental neighborhood system for $LAG \cap G(p,T\bo)\subset  LAG$.  The hypothesis of strict  \lag\ convexity for $\bo$ implies that 
there exists $\overline \epsilon >0$ so that
  $\tr_W \Hess \rho \geq \overline\epsilon$ for all \lag\  planes $W$  at points of $\bo$ with $|\cos\theta|<\epsilon$
  for some $\epsilon>0$.   Consequently, we have by equation (\KK.4)  that 
   $$
  \tr_W \Hess (-\log\delta) \ \geq \ {\overline \epsilon \over 2\delta}
  $$
    near $\bo$ for all  \lag\ planes $W$  with $|\cos\theta|<\epsilon$.
    
  Now choose $M>>0$ so that $\tr_W \Hess \rho \geq -M$ in a neighborhood of $\bo$ for all 
  $W\in LAG$.  Then, by (\KK.4) 
 $$
  \tr_W \Hess (-\log\, \delta) \ \geq\   -{M\over \d}+  {1\over \d^2}  |\nabla \rho\hk W|^2.
  $$
  If $|\cos\theta|\geq \epsilon$, this 
   is positive in a neighborhood of $\bo$ in $\Omega$.  This proves that $-\log \d$ is strictly \lagpsh\ near $\bo$.    By Theorem \JJ.6 the domain $\O$ is strictly $\f$-convex at infinity.
  \qed  
 
 \medskip
 
 Although a defining function for a strictly \lag\ convex boundary may not be \lag-plurisubharmonic, we do have the 
 following.
 
 \Prop{\KK.6} {\sl  Suppose $\O\subset\subset X$ has strictly \lag\ convex boundary $\bo$ 
 with defining function   $\rho$.  Then, for $A$ sufficiently large, the function 
 $\overline \rho \equiv \rho + A \rho^2$ is strictly Lag-psh
  in a neighborhood of $\bo$  and also a defining function for $\bo$.}

 \pf 
 By Lemma \KK.3  we have
 $$
 \tr_W (\Hess \overline\rho)\ =\ (1+2A\rho) \tr_W(\Hess \rho) +2 A |\nabla \rho\hk W|^2.
 \eqno{(\KK.5)}
 $$
  As noted in the proof 
 of Theorem \KK.5, there exist $\e, \overline\e>0$ so that $\tr_W (\Hess \rho) \geq \overline \e$ 
 for $W\in LAG$ with
 $|\cos\theta(W)|<\e$, by the strict boundary convexity.  Therefore $\tr_W(\Hess \overline\rho) \geq (1+2A\rho)\overline \e$ if $W\in LAG$ with $|\cos\theta(W)|<\e$.  Choose a lower bound $-M$ for $\tr_x (\Hess \rho)$ over all $W\in LAG$ for a neighborhood of $\bo$.  Then by (\KK.5), 
 $\tr_W (\Hess \rho) \geq -(1+2A\rho)M + 2|\nabla\rho|^2 A \e^2$ for $W\in LAG$ with 
 $|\cos \theta(W)|\geq \e$.  For $A$ sufficiently large, the right hand side is $>0$ in some neighborhood of $\bo$.\qed
 \medskip
 
 This leads to the following, which will be useful in the next section.
 
 \Theorem{\KK.7}  {\sl  Suppose $\O\subset\subset X$ has strictly \lag\ convex boundary and that 
 there exists a smooth strictly Lag-psh function $u$ on $\ob$.  Then $\O$ admits a strictly 
 Lag-psh defining function.}
 \medskip
 Note that in the case $\O\ss\ss X\equiv\bbc^n$ there many strictly Lag-psh functions on $X$.
  
  \pf
  By Proposition \KK.6 there exists a smooth \lag psh function $\rho$ defined on a neighborhood of $\bo$.
  For $\d_1 > > \d_2 >0$ sufficiently small, the function $\max\{\rho, -\d_1 +\d_2 u \}$ is Lag-psh
  and equal to $\rho$ near $\bo$.  We now apply the maximum smoothing (see pages 373-4 in [\PinGGC]) to obtain the desired defining function.\qed
 \medskip

  The (strict)  \lag\ convexity of a boundary can be equivalently defined in terms of its second fundamental form.  Note that if $Y\subset X$ is a smooth hypersurface with  a chosen unit normal field $n$ we have the associated {\bf second fundamental  form} $II$ defined on $TY$ by 
  $$
  II(v,w) \ \equiv \ \langle  \nabla_v \wt w, n \rangle
  $$  
  where $\wt w$ is any extension of $w$ to a tangent vector field on $Y$.
  For example, when $H=S^{n-1}(r)\subset \bbr^n$ is the euclidean sphere of radius $r$, oriented by the outward-pointing unit normal,  we find that $II(V,W) = -{1\over r} \langle V,W\rangle$.

Recall the following standard fact (cf. [\PTCG, Lemma 5.11]).

\Lemma{\KK.8} {\sl Suppose $\rho$ is a defining function for $\O$ and let II denote the second
fundamental form of the   hypersurface $\bo$ oriented by the outward-pointing normal.
Then 
$$
\Hess\, \rho\,\bigl|_{T\bo} \ =\ -|\nabla \rho| \, II
$$
and therefore
$$
\tr_W \Hess \rho\ =\ -|\nabla \rho |\, \tr_W II
$$
for all $W\in G(n, T\partial \O)$ and in particular for all  $W\in LAG\cap G(n, T\partial \O)$.}

  \Remark{}
    Recall that a defining function $\rho$ for $\Omega$ satisfies  $|\nabla \rho| \equiv 1$ in a neighborhood of $\bo$ if and only if $\rho$ is the signed distance to $\bo$ ($<0$ in $ \Omega$ and $>0$ outside of $\Omega$).   In fact any function $\rho$ with $|\nabla \rho| \equiv 1$  in a riemannian manifold is, up to an additive constant,   the distance function to (any) one of its level sets.
  In this case the lemma implies that
      $$
    \Hess \rho \ =\ \left( \matrix{0 & 0 \cr 0 & - II}       \right)
   \eqno{(\KK.6)} 
   $$
   where $II$ denotes the second fundamental form of the hypersurface $H=\{\rho =\rho(x)\}$
   with respect to the normal $n= \nabla \rho$ and the blocking in (\KK.6) is with respect to the 
   splitting $T_xX = \span (n_x) \oplus T_x H$.     
   
  As an immediate consequence of Lemma \KK.8 we have
  
  \Prop{\KK.9}  {\sl  Let $\Omega\subset X$ be a domain with smooth boundary $\bo$ oriented by
  the outward-pointing normal. Then $\bo$ is \lag\ convex if and only if its second fundamental form satisfies
  $$
  \tr_W II \ \leq \ 0 
    \eqno{(\KK.7)} 
   $$
  for all \lag\ planes $W$ which are tangent to $\bo$, i.e., for all $W\in LAG(\bo)$.  This can be expressed more geometrically by saying that
  $$
   \tr_W B \ \ \ {\rm must\ be\  inward-pointing }
  $$
  for all $W\in LAG(\bo)$, where $B_{v,w} \equiv (\nabla_v\wt w)^{Nor}$ is the normal-vector valued second fundamental form.
  
  Furthermore, $\bo$ is  strictly \lag\ convex if and only if 
  $$
  \tr_W II\ <\ 0 \qquad\forall\, W\in LAG(\bo),
 \eqno{(\KK.8)} 
   $$
   or equivalently, $\tr_W B$ is non-zero and inward-pointing for all $W\in LAG(\bo)$.
   }

  \Remark{\KK.10}
  If $\rho $ is the signed distance to $\bo$, then equation (\KK.6) together with Lemma \KK.8
  can be used to simplify (\KK.4).  An arbitrary $n$-plane $W$ at a point can be put in the canonical form $W=(\cos\theta n + \sin \theta e_1)\wedge e_2\wedge \cdots \wedge e_n$ with $n=\nabla \rho$
  and $n, e_1,...,e_n$ orthonormal.  Then $\eta=e_1\wedge \cdots \wedge e_p$ is the tangential
  projection of $W$.  Note that $\tr_W \Hess \rho = -\sin^2\theta \tr_\eta II$ and that 
  $|\nabla \rho\hk W|^2=\cos^2\theta$, so that (\KK.4) becomes
  $$
  \tr_W \Hess (-\log \rho) \ =\ -{1\over \delta}\sin^2\theta\, \tr_\eta II +{1\over \delta^2} \cos^2\theta
  $$

\def\hol{\ch}

 Let $n$ denote the outward unit normal field to $\bo$.
 Then at each point $x\in \bo$ we have an orthogonal decomposition
 $$
 T_xX\ =\ \bbr n\oplus T_x \bo\ =\ \bbr n\oplus \bbr Jn\oplus \hol
 $$
  where $\hol=T_x(\bo) \cap JT_x(\bo)$ is the (unique) maximal complex subspace  
 in $T_x(\bo)$.
  
  \Prop{\KK.11}  {\sl  Let $W\subset T_xX$ be any Lagrangian $n$-plane at a point $x\in \bo$.
  Then $W$ is of the form
  $$
  W\ =\ (n\cos \theta +Jn\cos \theta)\wedge W_0
  $$
  where $W_0$ represents a \lag\ $(n-1)$-plane in the complex subspace $\hol$.
  In particular, every tangential \lag\ $n$-plane is of the form}
  $$
   W\ =\ Jn\wedge W_0
  $$
  
  \pf   Put  $W$   in  canonical form $W=(\cos\theta n + \sin \theta e_1)\wedge e_2\wedge \cdots \wedge e_n$ as above  with $n, e_1,...,e_n$ orthonormal.  Set
  $e(\theta) = \cos\theta n + \sin \theta e_1$.  If $W$ is Lagrangian, then
  $e(\theta), Je(\theta)$,  $e_2, Je_2,...,e_n, Je_n$ form an orthonormal basis of $\bbr^{2n}$.
  In particular, $e_2, Je_2,...,e_n, Je_n$ are perpendicular to $\span\{e(\theta), Je(\theta)\}
  = \span\{n,e_1\}$. Hence, $\span\{n,e_1\}$ is a $J$-invariant.  We conclude that $e_1=\pm Jn$
  and $\hol=\span\{e_2, Je_2,...,e_n, Je_n\}$.\qed
  \medskip
  
  Combined with Proposition \KK.9 we conclude the following, which provides another way of describing Lagrangian boundary convexity.
  
  \Cor{\KK.12}  {\sl  Let $II$ be the  second fundamental form of $\bo$ with respect to the outer unit normal field $n$
as above.  Then $\bo$ is \lag\ convex if and only if 
  $$
  II(Jn, Jn)+\tr_{W_0} II \ \leq \ 0
 \eqno{(\KK.9)} 
 $$
  for all \lag\ $(n-1)$-planes $W_0$ in the holomorphic tangent space $\hol$.
  Furthermore,  $\bo$ is  strictly \lag\ convex if and only if  the inequality in (\KK.9) is strict for all $W_0$. }


\vskip .3in
\centerline{\bf  \LL.  The Dirichlet problem for the \lag\ Monge-Amp\'ere equation.}
\medskip

It is an important fact that Lagrangian harmonic functions exist in abundance locally on
any Gromov manifold $X$.  In fact the Dirichlet problem can always be solved on any domain
$\O\ss \ss X$ with smooth boundary such that
$$
\O \ \ {\sl admits \ a \ strictly \ Lagrangian \ plurisubharmonic \ defining \ function. }
\eqno{(\LL.1)}
$$
 By Theorem \KK.7 this is equivalent to the fact that 
 $$
 \eqalign
 {
&\bo \ \ {\sl is \ strictly\  \lag\ convex \  and \ there\ exists}  \cr
{\sl   a \  smo}  &{ \sl oth\ strictly \ \lag \ plurisubharmonic \ function \ on\ \  } \ob.
}
\eqno{(\LL.1)'}
$$
We make this assumption in each of the following theorems.
Note that all sufficiently small metric balls about any point have this property, 
and also that any domain $\O\ss\cn$ admits a smooth strictly Lag-psh function.

\Theorem {\LL.1.  (The Homogeneous  Dirichlet Problem)} {\sl
For any continuous $\vf \in C(\bo)$
there exists a unique function $u\in C(\ob)$ such that 
\medskip

(1) \ \ $u\bigr|_{\O}$ is \lag\  harmonic on $\O$, and
\medskip

(2)  \ \ $u\bigr|_{\bo} = \vf$.
}

\pf
This is a special case of Theorem 16.1 with $\GG=\LAG$  in [\DDR]. 
(It  follows as well from the  more general Theorem 13.1 in [\DDR].) \qed
\medskip

We can also treat the inhomogeneous equation  for the Lagrangian Monge-Amp\'ere operator:
$$
\mlag(\Hess \, u) \ = \  \psi,
\qquad
u\ \ {\rm Lag-psh}
$$
with continuous inhomogeneous term $\psi$.  Existence and uniqueness are easier when $\psi >0$  and smooth, and we outline below the proof given in [\DDR].
The case where $\psi\geq 0$, that is, where $\Hess \, u$ is allowed to hit the boundary of $\cp(\LAG)$,
and $\psi$ is continuous, is more complicated, and we shall refer to [\IDP] for  that case.

\Theorem {\LL.2.  (The  Inhomogeneous Dirichlet Problem)} {\sl
 Fix a continuous $\psi \geq 0$ on $\ob$.
 Then for any  $\vf \in C(\bo)$
there exists a unique function $u\in C(\ob)$, which is Lagrangian plurisubharmonic on $\O$,   such that 
\medskip

(1) \ \ $u\bigr|_{\O}$ is the viscosity solution of $\mlag(\Hess \, u) = \psi$ on $\O$, and
\medskip

(2)  \ \ $u\bigr|_{\bo} = \vf$.
}

\pf
We shall first consider smooth $\psi > 0$. Note to begin that  
$$
\{A\in \cp(\LAG) : \mlag(A)\geq 1\}
$$
is unitarily invariant.  As a result this ``universal''  subequation determines the subequation
$\mlag(\Hess\, u ) \geq 1$ on any Gromov manifold
(see Chapter  5  in [\DDR].)
Theorem 13.1 in [\DDR] then applies to give the version of Theorem \LL.1 for the equation
$\mlag(\Hess\, u ) = 1$.   

Now  by conjugating $\Sym(T^*X)$ by 
$$
A\ \  \mapsto \ \ \left( \psi^{-2^{n/2}}\right) A \left(  \psi^{-2^{n/2}} \right) ^t  =  \psi^{-2^{n} } A,
$$
we get a jet-equivalence of  our equation $\mlag(\Hess\, u ) = \psi$
with the  equation  $\mlag(\Hess\, u ) = 1$.
We can then apply Theorem 13.1$'$ in [\DDR].  

For continuous $\psi\geq 0$ the reader is referred to [\IDP, Thm. 8.1].   \qed

\medskip

We now turn attention to the branches  $\cp_k(\LAG)$  of $\mlag$ which determine subequations
on any Gromov manifold $X$ (see Remark \HH.4 and Note \CC.9).

\Theorem {\LL.3.  (The Homogeneous Dirichlet Problem for the Branches)} {\sl
Fix a  continuous function $\vf \in C(\bo)$
and any $k=1, ... , 2^n$.  Then 
there exists a unique function $u\in C(\ob)$ such that 
\medskip

(1) \ \ $u\bigr|_{\O}$ is $\cp_k(\LAG)$-harmonic on $\O$, and
\medskip

(2)  \ \ $u\bigr|_{\bo} = \vf$.
}

\pf
This is Theorem 13.1 in [\DDR].  We point our that $\cp_k(\LAG)$ is a riemannian U$(n)$-subequation
on $X$.  Furthermore, the boundary is strictly $\cp(\LAG)$-convex, which implies strict $\cp_k(\LAG)$-convexity
and strict  dual $\cp_k(\LAG)$-convexity (= strict  $\cp_{2^n-k}(\LAG)$-convexity) since $\cp(\LAG)
=\cp_1(\LAG) \ss\cp_k(\LAG)$ for all $k$.  So the boundary condition is satisfied.  \qed

\Note{\LL.4}  Theorem \LL.3 remains true under the  weaker natural boundary convexity condition,
(discussed in the proof above)  that: \smallskip
 
 \centerline
 {\sl There exists a strictly $\cp_{\ell}(\LAG)$-subharmonic defining function for $\bo$}

 \centerline{
where $\ell = \min\{k, 2^n-k\}$.
}
  This is equivalent (as above for $k=1$)
 to $\bo$ being strictly $\cp_{\ell}(\LAG)$-convex
and the existence of a strictly $\cp_{\ell}(\LAG)$-subharmonic  function on $\ob$.

\Note{\LL.5}   There is also an operator ${\rm M}_{\Lag,k} \equiv \L_k \L_{k+1} \cdots \L_{2^n}$ on
the subequation $\cp_k(\LAG)$,  where $\L_1 \leq \L_2 \leq \cdots$ are the ordered 
eigenvalues of $\mlag= {\rm M}_{\Lag,1}$ (cf. \S \CC),
although it is not a polynomial operator unless $k=1$.
This operator is discussed,
for example,  in [\IDP]; see in particular,  (6.3)  and Prop. 6.11 in [\IDP].
The inhomogeneous Dirichlet Problem for ${\rm M}_{\Lag,k}$  
can be solved  uniquely for any inhomogeneous term $\psi\geq 0$
under the boundary assumptions in Note \LL.4, again by the Theorem 8.1 in [\IDP].


\vskip .3in
\centerline{\bf  \MM.  Ellipticity of the Linearization.}
\medskip

It is natural to consider the linearization of the operator  ${\rm M}_\Lag(\Hess\, u)$ on compact subsets
of the interior of its domain, i.e., on compact subsets of $\Int \cp(\LAG)$.  
Of course for any weakly elliptic operator $f$ (one where $f(A+P)\geq f(A)$ for $P\geq0$ and $A$ in its domain),
the linearization $L(f)$ is  weakly elliptic;  and if $f$ is uniformly elliptic, so is $L(f)$. 
However, there is a  (not particular trivial) result  that for operators defined by
G\aa rding/Dirichlet  polynomials, such as $\mlag$, the linearization at all interior points
of the subequation is always positive definite.  Details of this can be found in 
[\IDP] and [\HLGG].

\Prop{\MM.2}  {\sl Let $\O\ss X$ be a compact domain, and assume that 
$u$ is a $C^2$ $\LAG$-plurisubharmonic function on $\ob$ which satisfies the equation 
$$
\mlag (\Hess u)\ =\ f\ >\ 0
$$
on $\ob$.  
Then the linearization of $\mlag $ at this solution is 
uniformly elliptic.
}
\medskip
This  allows one to use Implicit Function Techniques to get smooth solutions for nearby boundary
data.

\vskip.3in


\centerline {\bf Appendix A. }\smallskip
\centerline{ \bf A More Detailed Presentation of the Lagrangian Subequation.}   
\medskip

A more detailed algebraic discussion of the subequation $\cp(\LAG)$ is presented
here. It includes a description of the extreme rays. The material in this appendix is self-contained
and includes a second treatment of the results of Section \BB.

\def\CCH{{\rm CCH}}

First, for the reader's convenience, we summarize results that hold for any geometric
subequation (see ... for more details).  Start with any closed subset $\GG$
(say $\GG=\LAG$) of the Grassmannian $G(p, \bbr^N)$ of unoriented $p$-planes in $\bbr^N$.
Identify $W\in \GG$ with orthogonal projection $P_W\in\Sym(\bbr^N)$ onto $W$.
We have the following concepts.

\medskip
\noindent
(A.1) {\bf (The Subequation $\cp^+$)}\ \ $A\in\cp^+ \ \ \iff\ \ \tr\left(A\bigr|_W\right) = \bra A {P_W} \geq 0 \ \ 
\forall\, W\in\GG$.

\medskip
\noindent
(A.2) {\bf (The Convex-Cone Hull $\cp_+$)}\ \ $\cp_+  \equiv   $ the convex-cone hull of $\GG$,
denoted $\CCH(\GG)$.

\medskip
\noindent
(A.3) {\bf (Polars)} \ \   $\cp^+$ and $\cp_+$  are polar cones in $\Sym(\bbr^N)$.

\medskip
\noindent
(A.4) {\bf (The Edge $E$)} \ \ $E\ \equiv \ \cp^+ \cap (-\cp^+)$,  that is,   

\qquad\qquad\qquad
$A\in E \iff \tr\left (A\bigr|_W\right) = \bra A {P_W}=0 \ \ 
\forall\, W\in\GG$.

\medskip
\noindent
(A.5) {\bf (The Span $S$)}\ \ $S\ \equiv \ \span \GG \ \equiv\  \span \cp_+$.

\medskip
\noindent
(A.6) {\bf ($E=S^\perp$)} \ \ $\Sym(\bbr^N) = E\oplus S$ is an orthogonal decomposition,

\medskip
\noindent
(A.7)  {\bf (Extreme Rays in $\cp_+ \equiv \CCH(\GG)$)} \ \ The extreme rays in $\cp_+$  are the rays through

\quad $P_W$ where $W\in \GG$.
\medskip

Note that (A.1), (A.2)  and (A.5) are definitions,
while (A.4) is a definition combined with  an immediate consequence of (A.1).
The polar facts (A.3) and(A.6) are easy. For (A.7) note that from the definition of $\cp_+$
every extreme ray is generated by a $P_W$. On the other hand, all of the elements
$P_W$ lie on a sphere  in the hyperplane $\{\tr A = p\}$, centered 
about  ${p\over n}{\rm I}$, and it then follows that every $P_W$ is extreme.

\medskip
The more interesting geometric cases tend to have a non-trivial edge.
The edge $E$ is a vector subspace of $\cp(\GG)$ and it contains all 
other vector subspaces of $\cp(\GG)$.  It can be ignored in the following sense.
Let $\pi: \Sym(\bbr^N) \to S$ denote orthogonal projection.

\medskip
\noindent
(A.8a)  {\bf  (The Reduced Subequation $\cp_0^\perp$)}  \ \ $\cp^+ = E\oplus \cp_0^+$ defines 
$ \cp_0^+$, and

\medskip
\noindent
(A.8b) \ \  $\cp_0^+ = \pi(\cp^+) = \cp^+\cap S = \cp^0_+$ \ (= the polar of $\cp_+$ in $S$).

\bigskip

Now we give an explicit description of each of the objects above in the case at hand, namely
$\GG=\LAG$. We use the standard fact that the action of the unitary group on 
$  [{\rm I}]  \oplus \Herm^\sk(\cn)$ has a cross-section $D$ with each orbit intersecting
$D$ in a finite number of points.

\medskip
\noindent
(A.9) \ \ $D \equiv \left \{H(t,\l) \equiv {t\over 2n}{\rm I} + \sum_j \l_j(P_{e_j} -P_{Je_j}) : (t,\l)\in \bbr^{n+1}\right\}$.

\medskip
\noindent
(The label $D$ is used since as a $2n\times 2n$-matrix, each element of $D$ is diagonal.)
\smallskip

The {\bf Lagrangian part} of $A\in \Sym_\bbr(\cn)$ is defined to be

\medskip
\noindent
(A.10a) \ \ $A^\LAG\ \equiv\ {1\over 2n} (\tr A){\rm I} + \half(A+JAJ)\ \equiv \pi(A)$,
\medskip

and the {\bf skew-hermitian part} of $A$ is
\medskip
\noindent
(A.10b) \ \ $A^\sk \ \equiv \ \half(A+JAJ)$,
\medskip

Let $W(\ve) = W(\pm  1, ... ,\pm 1) \equiv W(\ve_1, ... ,\ve_n)$ be the  axis Lagrangian $n$-plane
defined by the condition that  $e_j \in W(\ve)$ if $\ve_j=+1$ and $Je_j\in W(\ve)$ if $\ve_j=-1$. Then

\medskip
\noindent
(A.11)\ \ $P_{W(\ve)} = \half{\rm I} + \half\sum_j\ve_j (P_{e_j} -P_{Je_j}) $ 
\medskip
 \noindent
 In particular, these $2^n$ axis Lagrangian planes $P_{W(\ve)} \in D$ (with $t=n$ and $\l={\ve\over 2}$)
comprise the vertices of the cube $[-\half, \half]^n$ in the $t=n$ hyperplane.
\medskip

 It is easy to see that $\GG\cap D = \LAG\cap D$ consists of the $2^n$ points $P_{W(\ve)}$. 
This proves 

\medskip
\noindent
(A.12a)\ \ $\span G\cap D = D$,  \  and hence

\medskip
\noindent
(A.12b)\ \ $S \ \equiv \ \span \GG =   [{\rm I}] \oplus \Herm^\sk(\cn)$.

\medskip

Thus, all of the facts in (A.1) through (A.8) are true with the edge $E$ of the subequation $\cp^+$
identified as 

\medskip
\noindent
(A.13)\ \ $E\ =\ \Herm_0^\sym(\cn)$, the space of traceless hermitian-symmetric forms, 
\medskip
\noindent
since this space equals
$S^\perp = (  [{\rm I}] \oplus \Herm^\sk(\cn))^\perp$.
\medskip

Since the $2^n$ points $P_{W(\ve)}$ are the vertices of a cube in the affine hyperplane $\{t=n\}\cap D$ in $D$,
the convex hull, which is $\cp_+\cap D$, equals the cone on this cube.
Consequently, by (A.11),

\medskip
\noindent
(A.14a)\ \ $\cp_+\cap D = \{H(t,\l) : \sup_j |\l_j|\leq {t\over 2n}\}$ with $P_{W(\ve)} = H(n, {\ve\over 2})$, \ so that
\medskip
\noindent
(A.14b)\ \   $P_{W(\ve)} \in \partial (\cp_+\cap D)$ generate the extreme rays.
\medskip

The polar of this cone $\cp_+$ in $D$, which equals $\cp^+\cap D$, is now easy to compute.
First compute that 

\medskip
\noindent
(A.15)\ \ $\bra{H(t,\l)}{P_{W(\ve)}} \ =\ {t\over 2} + \sum_j \ve_j \l_j$.  
\medskip
\noindent
This term  is $\geq 0 \ \forall \ \ve=(\pm 1, ... ,\pm 1)\ \iff \ |\l_1|+\cdots+|\l_n|\leq {t\over 2}$.  Therefore,

\medskip
\noindent
(A.16a)\ \ $\cp^+\cap D \ =\ \{ H(t,\l) :  |\l_1|+\cdots+|\l_n|\leq {t\over 2}\}.$

\medskip
This proves that the subequation $\cp^+$ is given by

\medskip
\noindent
(A.16b)\ \  $\cp^+ = \{A\in \Sym_\bbr(\cn) : {\tr A\over 2} - \l_1^+ - \cdots - \_n^+\geq 0\}$, 
\medskip\noindent
where
$\l_1^+,...,\l_n^+$ are the non-negative eigenvalues of $A^\sk$.

\medskip
Setting $t=1$ in (A.16a), we see that $\cp^+\cap D$ is the cone on the set $|\l_1|+\cdots + |\l_n|\leq\half$.
Hence, its extreme rays are through $H(1,{\a\over 2})$ where $\a$ is one of the $2n$ unit vectors on an
axis  line in $\rn$.  That is:

\medskip
\noindent
(A.17a)\ \  The elements $H(1,{\a\over 2}) \equiv {1 \over 2n} {\rm I} \pm \half (P_{e_j}- P_{J e_j})$ (for some $\pm$ and $j=1,...,n$),

\hskip .3in generate the extreme rays in $\cp^+\cap D$.
\medskip
\noindent
On the other hand,
\medskip
\noindent
(A.17b)\ \ $P_{e_j}^\LAG \ =\ {1\over 2n}{\rm I} +\half(P_{e_j}- P_{J e_j})$ and 
$P_{Je_j}^\LAG \ =\ {1\over 2n}{\rm I} -\half(P_{e_j}- P_{J e_j})$.
\medskip

This proves the following.

\Theorem{A.1} {\sl
The extreme rays in the reduced subequation $\cp_0^+$ are through $P_e^\LAG$, $e\neq 0$
}
\medskip

This has important special consequences for $\GG\equiv \LAG$.

\Prop{A.2}  {\sl
Let $\GG \equiv \LAG$.  Then, regarding $\cp^+$, one has:
\medskip

(1)\ \ $\cp^+_0\ =\ \pi(\cp), \quad \Int\cp^+_0 =\pi(\Int \cp)$,

\medskip

(2)\ \ $\cp^+ \ =\ E+\cp \ \ = \Herm^\sym_0(\cn)+\cp$,

\medskip

(3)\ \ $\Int \cp^+ \ =\ E + \Int \cp \ = \Herm_0^\sym(\cn)+\Int\cp$.

\medskip
\noindent
While for $\cp_+$ one has

\medskip

(4)\ \ $\cp_+ \ =\ \cp\cap S \ = \cp \cap (  [{\rm I}]  \oplus \Herm^\sk(\cn), \and
          \Int_{\rm rel} \cp_+ = (\Int \cp) \cap S.$

}

\pf
Assertion (1)  is immediate, and assertions (2) and (3) follow.
Since $(\cp\cap S)^0 = \cp^0 + S^0 = \cp+E$, which equals $\cp^+$ by (2), 
this proves (4). For a second, more direct proof of (4) note that $H(t,\l)\geq0 \ \iff\ |\l_j|\leq {t\over 2}, j=1,..,n
\ \iff\ H(t,\l) \in \cp_+\cap D$ by the definition of $H(t,\l)$ and (A.14) respectively.
\qed

\Remark {A.3} By contrast to Proposition A.2, in the cases 
$\GG= G_K(1, K^n)\ss\Sym_\bbr(K^n)$ for $K=\bbr, \bbc$ or  $\bbh$,
the convex-cone hull $\cp_+  \equiv \CCH(\GG)$ and the reduced subequation
$\cp^+_0 = \cp^+\cap S$ are self-polar in $S\equiv \span \GG$.
 \medskip

Note that the edge $E$ generates $\cp^+$ as a subequations, i.e., 
\medskip

 (2)\ \ $\cp^+ = E+\cp$ in all four  cases $\GG=\LAG$ and $\GG=G_K(1, K^n)$ for $K=\bbr, \bbc, \bbh$.

\Remark{A.4}  For any subspace $E\ss \Symn$ with $E\cap \cp = \{0\}$,
the sum $E+\cp$ is closed and hence a subequation.  
In fact, it has edge $E$, satisfies Theorem A.1, and enjoys all the
properties in Proposition A.2.  All of this is proven in [\EDGE].

\vskip.3in

\centerline{\bf References}
\medskip

\item{[BT$_1$]}  E. Bedford and B. A. Taylor,  
 The Dirichlet problem for a complex Monge-Amp`ere equa- tion, Inventiones Math.37 (1976), no.1, 1-44.
 \smallskip
 
 \item{[BT$_2$]} \ -------- \ , \ 
Variational properties of the complex Monge-Amp`ere equation, I. Dirichlet principle, Duke Math. J. 45 (1978), no. 2, 375-403.
\smallskip

 \item{[G\aa]}  L. G\aa rding,  {\sl An inequality for hyperbolic polynomials}, J. Math. Mech. {\bf 8}  no. 2 (1959), 957-965.

\smallskip

\item{[G]}  M. Gromov, {\sl Pseudoholomorphic curves in symplectic manifolds}, Invent. Math. {\bf 82} no. 2 (1985), 307-347.
\smallskip

\item {[\PTCG]}  F. R. Harvey and H. B. Lawson, Jr.,  {\sl  An introduction to potential theory in calibrated geometry},    Amer. J. Math.  {\bf 131} no. 4 (2009), 893-944.  ArXiv:math.0710.3920.

\smallskip

\item {[\DDD]} \  -------- \ , \  {\sl  Dirichlet duality and the non-linear Dirichlet problem},
    Comm. on Pure and Applied Math. {\bf 62} (2009), 396-443.

\smallskip

\item {[\PinGGC]}   \  -------- \ , \ {\sl  Plurisubharmonicity in a general geometric context},
  Geometry and Analysis {\bf 1} (2010), 363-401. ArXiv:0804.1316

\smallskip

\item {[\DDR]}  \  -------- \ , \  {\sl   Dirichlet duality and the non-linear Dirichlet problem on Riemannian manifolds},
  J. Diff. Geom.  {\bf 88} No. 3 (2011), 395-482.  ArXiv:0907.1981.

\smallskip

\item {[\HLG]}  \  -------- \ , \ {\sl  Hyperbolic polynomials and the Dirichlet problem},
    ArXiv:0912.5220.

\smallskip

\item {[\GPandC]} \  -------- \ , \  {\sl  Geometric plurisubharmonicity and convexity - an introduction},
  Advances in Math.  {\bf 230} (2012), 2428-2456.   ArXiv:1111.3875.

\smallskip

\item {[\HLGG]}   \  -------- \ , \  {\sl  G\aa rding's theory of hyperbolic polynomials},
   {\sl Communications in Pure and Applied Mathematics}  {\bf 66} no. 7 (2013), 1102-1128.



\smallskip

\item {[\SURVEY]}  \  -------- \ , \  {\sl  Existence, uniqueness and removable singularities
for nonlinear partial differential equations in geometry},
  pp. 102-156 in ``Surveys in Differential Geometry 2013'', vol. 18,  
H.-D. Cao and S.-T. Yau eds., International Press, Somerville, MA, 2013.
ArXiv:1303.1117.

\smallskip

\item {[\REST]}  \  -------- \ , \  {\sl  The restriction theorem for fully nonlinear subequations},
     Ann. Inst.  Fourier, {\bf 64}  no. 1 (2014), p. 217-265.  ArXiv:1101.4850.

\smallskip

\item {[\RSing]}  \  -------- \ , \  {\sl  Removable singularities for nonlinear subequations},
   Indiana Univ. Math. J., {\bf 63}, No. 5 (2014), 1525-1552.
 ArXiv:1303.0437.

\smallskip

\item {[\ACM]}  \  -------- \ , \  {\sl  Potential theory on almost complex manifolds},
     Ann. Inst.  Fourier,  Vol. 65 no. 1 (2015), p. 171-210.  ArXiv:1107.2584.

\smallskip

\item  {[\AE]}  \  -------- \ , \ {\sl The AE Theorem and Addition Theorems for quasi-convex functions,}  \ ArXiv:1309:1770.

\smallskip

\item {[\IDP]}   \  -------- \ , \  {\sl  The inhomogeneous Dirichlet Problem},
  (to appear).

\smallskip

\item  {[\EDGE]}  \  -------- \ , \ {\sl  Pluriharmonics in general potential theories,}  (to appear).

\smallskip

\item  {[\CONT]}  \  -------- \ , \ {\sl  Contact sets for general subequations,}  (to appear).

\smallskip

\item{[HW$_1$]}  F. R. Harvey and R. O. Wells, Jr.,  {\sl Holomorphic approximation and hyperfunction theory on a $C^1$ totally real submanifold of a complex manifold},  Math. Ann. {\bf 197} (1972), 287-318. 

\smallskip
\item{[HW$_2$]}  \  -------- \ , \  {\sl Zero sets of non-negative strictly plurisubharmonic functions},   Math. Ann. 
{\bf 201} (1973), 165-170.

\smallskip

\item{[LM]} H. B. Lawson and M.-L. Michelsohn, Spin Geometry,  Princeton Math. Series, 38. Princeton Univ. Press, Princeton, NJ, 1989.

\smallskip

\item{[MS]}  D. McDuff and D. Salamon,  J-holomorphic curves and symplectic topology. American Mathematical Society Colloquium Publications, 52. American Mathematical Society, Providence, RI, 2012.

 \smallskip

 \item{[M]} J. Milnor, Morse Theory, Annals of Math. Studies no. 51, Princeton Univ. Press, Princeton, NJ,  1963.
 
\item {[NW]}  A. Nijenhuis and W. Woolf, {\sl  Some integration problems in almost -complex and complex manifolds}, 
Ann. of Math.,  {\bf 77} (1963), 424-489.

\smallskip

\item {[O]}  A. Oberman, {\sl The  convex envelope is the solution of a nonlinear obstacle problem}, 
Proc. A.M.S. {\bf 135}  (2007), no. 6,  1689-1694.

\smallskip

\item {[OS]}  A. Oberman and L. Silvestre, {\sl The Dirichlet problem for the convex envelope}, 
Trans. A.M.S. {\bf 363}  (2011), no. 11,  5871-5886.

\smallskip

\end